\documentclass[smallextended,numbook,runningheads]{svjour3}
\smartqed
\usepackage{mathptmx}
\usepackage{amssymb}
\usepackage{mathrsfs}
\usepackage{amsfonts,dsfont}
\usepackage{graphicx}
\usepackage{amsmath}
\usepackage{psfrag}
\usepackage{booktabs}
\usepackage{cite}
\usepackage{url}
\usepackage[usenames]{color}
\numberwithin{equation}{section}
\usepackage[top=1.3in, bottom=1.3in, left=1.2in, right=1.5in]{geometry}

\newcommand{\R}{\mathbb{R}}
\newcommand{\N}{\mathbb{N}}

\newcommand{\dd}{\text{d}}

\newtheorem{thm}{Theorem}[section]

\newtheorem{prop}[thm]{Proposition}

\newtheorem{assumption}[thm]{Assumption}

\journalname{Journal of Scientific Computing}
\begin{document}

\title{Optimal error estimates of Galerkin finite element methods for stochastic Allen-Cahn equation with additive noise
\thanks
{
RQ was supported by  NSF of China (11701073).
XW was supported by NSF of China (11671405, 11571373, 91630312),
NSF of Hunan Province  (2016JJ3137), Innovation Program of Central South University (2017CX017)
and Program of Shenghua Yuying at CSU. 
The authors want to thank the Tianyuan Mathematical Center in Northeast China
for the hospitality and Prof. Xiaobing Feng from The University of Tennessee for his useful comments
when this work was presented in a conference in June of 2018, hosted by the center.
}
}

\author{
Ruisheng Qi,
\,
Xiaojie Wang
}
\institute{
Ruisheng Qi \at School of Mathematics and Statistics, Northeastern University at Qinhuangdao, Qinhuangdao, China\\
 \email{qirsh@neuq.edu.cn} \\
 Xiaojie Wang \at School of Mathematics and Statistics, Central South University, Changsha, China\\
 \email{x.j.wang7@csu.edu.cn, x.j.wang7@gmail.com}
}


\date{Received: date / Accepted: date}

       \maketitle

       \begin{abstract}
          Strong approximation errors of both finite element semi-discretization and spatio-temporal full discretization
are analyzed for the stochastic Allen-Cahn equation driven by {\color{black}{additive  noise}} in space dimension $d \le 3$.
The full discretization is realized by combining the standard finite element method
with the backward Euler time-stepping scheme.
Distinct from the globally Lipschitz setting,  the error analysis becomes rather challenging and demanding,
due to the presence of the cubic nonlinearity in the underlying model.
By introducing two auxiliary approximation processes,  we
{\color{black}{propose an}} appropriate decomposition of the considered error terms
and introduce a novel approach of error analysis, to successfully recover
the convergence rates of the numerical schemes.
The approach is {\color{black}{original}} and does not rely on high-order spatial regularity properties of the approximation processes.
It is shown that the {\color{black}{fully}} discrete scheme possesses  convergence rates of order $ O(h^{\gamma} ) $ in space and
order $ O( \tau^{ \frac{\gamma}{2} }  )$ in time, subject to the spatial correlation of the noise process,
characterized by $ \|A^{\frac{\gamma-1}2}Q^{\frac12}\|_{\mathcal{L}_2}<\infty, \, \gamma {\color{black}{\in[\frac d3,2] }} $, ${\color{black}{d\in\{1,2,3\}}}$.
In particular,  a classical convergence rate of order $O(h^2 +\tau)$ is reachable, even in multiple space dimensions,
when the aforementioned condition is fulfilled with $ \gamma = 2 $.
Numerical examples confirm the previous findings.



\end{abstract}

\section{Introduction}

Stochastic partial differential equations (SPDEs)  are widely used to mathematically model random phenomena
appearing in the fields of physics, chemistry, biology, finance and many other branches of science.
Over the past decades, there have been plenty of research articles analyzing numerical discretizations of
parabolic SPDEs, see, e.g., monographs \cite{lord2014introduction,kruse2014strong}
and references therein. In contrast to an overwhelming majority of literature
focusing on numerical analysis of SPDEs with globally Lipschitz nonlinearity,
 only a limited number of papers investigated numerical SPDEs in the non-globally Lipschitz regime
\cite{becker2017strong,brehier2018strong,brehier2018analysis,brehier2018weak,
feng2017finite,gyongy2016convergence,jentzen2015strong,kovacs2015discretisation,
becker2016strong,hutzenthaler2016strong,jentzen2016exponential,liu2018strong,kovacs2015backward} and it is still far from {\color{black}{being}} well-understood.
As a typical example of parabolic SPDEs with non-globally Lipschitz nonlinearity,  stochastic Allen-Cahn equations,
perturbed by additive or multiplicative noises, have received increasing attention in the last few years.
Recently, several research works were reported on numerical approximations of such equations
\cite{feng2017finite,katsoulakis2011noise,kovacs2015backward,liu2017wong,
liu2018strong,brehier2018strong,Majee2017optimal,brehier2018weak,wang2018efficient,brehier2018analysis,kovacs2015discretisation}.
The present work makes further contributions in this direction,
by successfully recovering optimal strong convergence rates for finite element semi-discretization and spatio-temporal
full discretization of stochastic Allen-Cahn equations with additive noise, including both the space-time white noise in 
space dimension  $d = 1$ and the trace-class noise in multiple space dimensions.


%

Let $\mathcal{D} \subset \R^d, d \in \{ 1, 2, 3 \} $ be a bounded open spatial domain with smooth boundary
and let $H : = L_2 ( \mathcal{D}, \R )$ be the real separable Hilbert space endowed with usual inner product
and  norm.
Throughout  this article we are interested in the following semi-linear parabolic SPDE in $H$,
\begin{align}\label{SACE}
\begin{split}
\left\{\begin{array}{ll}
\dd X(t) + AX(t) \dd t
=F(X(t)) \dd t + \dd W(t),& (t,x)\in (0, T]\times \mathcal{D},
\\
X(0)=X_0,&x\in \mathcal{D},
\end{array}\right.
\end{split}
\end{align}
where $A \colon D ( A ) \subset H \rightarrow H$ is a linear, densely defined,
positive self-adjoint unbounded operator with compact inverse (e.g., $A = - \Delta$ with homogeneous Dirichlet boundary condition)
in $H$, generating an analytic semigroup $\mathcal{S}(t) = e^{ - t A }$ in $H$.
Moreover, $\{W(t)\}_{t\geq 0}$ is an $H$-valued (possibly cylindrical) $Q$-Wiener process on a filtered probability space
$(\Omega, \mathcal{F}, \mathbb{P}, \{\mathcal{F}_t\}_{t\geq 0})$ with respect to the normal filtration $\{\mathcal{F}_t\}_{t\geq 0}$.
The nonlinear mapping $F$  is assumed to be a Nemytskij operator, given by
$
F( u )( x ) = f( u ( x ) ),
\:
x \in \mathcal{D},
\text{ with }
f( v ) = v - v^3, v \in \R.
$
%
%
Such problem is often referred to as stochastic Allen-Cahn equation.
%
%
Under further assumptions specified later, particularly including
\begin{align}
\label{eq:intro-AQ-condition}
\|A^{\frac{\gamma-1}2}Q^{\frac12}\|_{\mathcal{L}_2}<\infty,\; \text { for  some }
\;
\gamma\in{\color{black}{[\tfrac{d}3,2]}},
\end{align}
it is proved in Theorems \ref{thm:uniqueness-mild-solution}, \ref{them:regulairty-mild-solution} that the problem \eqref{SACE} possesses a unique mild solution,
\begin{align}\label{eq:intro-mild-solution}
X(t)
=
\mathcal{S}(t)X_0
+
\int_0^t\mathcal{S}(t-s)F(X(s))\,\mathrm{d} s
+
\int_0^t\mathcal{S}(t-s) \, \mathrm{d} W(s),
\quad
t \in [0, T],
\end{align}
which enjoys the Sobolev and H\"{o}lder regularity properties
\begin{align}\label{eq:spatial-regulairty-mild-intr}
X \in L_\infty([0,T];L^{2p}(\Omega; \dot{H}^\gamma)),
\quad
\forall p\geq 1,
\end{align}
and
for $\forall p\geq 1$ and $ 0 \leq s < t \leq T $,
\begin{align}\label{eq:temporal-regulairty-mild-intr}
\|X(t)-X(s)\|_{L^{2p}(\Omega; \dot{H}^\beta)}
\leq
C(t-s)^{\frac{\min\{1,\gamma-\beta\}}2},
\quad
\beta\in[0,\gamma].
\end{align}
Here $\| \cdot \|_{\mathcal{L}_2} $ stands for the Hilbert-Schmidt norm, $\dot{H}^\alpha : = D(A^{\frac \alpha2}), \alpha\in \mathbb{R}$
and  the parameter $\gamma\in{\color{black}{[\frac{d}3,2]}}$ coming from \eqref{eq:intro-AQ-condition}
quantifies
the spatial regularity of the covariance operator $Q$ of the driving noise process
(Assumption \ref{assum:eq-noise}). {\color{black}{The setting covers both the space-time white noise in space dimension  $d = 1$ 
and the trace-class noise in multiple space dimensions (see Remark \ref{rem:sectoin2-noise} for details).}} 
The obtained space-time regularity coincides with that in \cite{kruse2012optimal} for SPDEs with globally Lipschitz nonlinearity
and is thus optimal in the spirit of \cite{kruse2012optimal}.
%

Let $\mathcal{D} \subset \R^d, d \in \{ 1, 2, 3 \} $, be an open convex polynomial domain and
$A=-\Delta$ with $D(A)=H^2 (\mathcal{D}) \cap H_0^1(\mathcal{D})$. Let $V_h$ be a finite element space
of piecewise continuous linear functions
and $X_h $ the finite element spatial approximation of the mild solution  $X$, which can be represented by
\begin{equation} \label{eq:intro-FEM-mild-form}
X_h ( t ) = \mathcal{S}_h(t)P_hX_0
+
\int_0^t\mathcal{S}_h(t-s)P_hF(X_h(s))\,\dd s
+
\int_0^t\mathcal{S}_h(t-s)P_h \,\dd W(s),
\quad
t \in [0, T].
\end{equation}
Here $\mathcal{S}_h(t):= e^{ - t A_h } $ is the strongly continuous semigroup generated by the discrete Laplace operator $A_h$.
The resulting spatial approximation error is measured  as follows (Theorem \ref{them-error-semi-discrete-problem})
\begin{align}\label{eq:error-semi-discrete-problem}
\|X(t)-X_h(t)\|_{L^{2p}(\Omega;H)}
=
O(h^\gamma),
\quad
\gamma\in {\color{black}{[\tfrac{d}3,2]}},
\end{align}
where $\gamma\in {\color{black}{[\frac{d}3,2]}}$ is determined by the assumption \eqref{eq:intro-AQ-condition}.
The obtained convergence rate in space is {\color{black}{called  optimal}} since it exactly coincides with the order of
optimal spatial regularity of the solution \cite[Chapter 1]{thomee2006galerkin}.
Discretizing the semi-discrete problem by a backward Euler time-stepping scheme, we also investigate
a {\color{black}{fully  discrete scheme}} for \eqref{eq:intro-mild-solution}, given by
\begin{align}
\label{eq:intro-one-step-backward-Euler}
\begin{split}
X_{h,m} = \mathcal{S}_{\tau,h} X_{h,m-1} + \tau \mathcal{S}_{\tau,h} P_hF(X_{h,m}) + \mathcal{S}_{\tau,h} P_h \Delta W_m,
\;
X_{h,0}=P_hX_0,
\;
m\in\{1,2,\cdots, M\},
\end{split}
\end{align}
where $X_{h,m}$ is the {\color{black}{fully}} discrete approximations of $X(t_m)$ and $\mathcal{S}_{\tau,h}:=( I + \tau A_h)^{-1}$.
{\color{black}{Equivalently, the one-step recursion \eqref{eq:intro-one-step-backward-Euler} can be reformulated as
\begin{equation}
X_{h,m}
=
\mathcal{S}^m_{\tau,h}X_{h,0}
+
\tau\sum_{i=0}^{m-1}\mathcal{S}^{m-i}_{\tau,h}P_hF(X_{h,i+1})
+
\sum_{i=0}^{m-1}\mathcal{S}^{m-i}_{\tau,h}P_h \Delta W_{i+1}.
\end{equation}}}
As stated in Theorem \ref{them:error-estimate-full-FME},
the corresponding strong approximation error reads,
\begin{align}
\|X(t_m)-X_{h,m}\|_{L^{2p}(\Omega;H)}
=O(h^\gamma
+\tau^{\frac\gamma2})
,
\quad
\gamma\in {\color{black}{[\tfrac{d}{3},2]}}.
\end{align}
This indicates how the strong convergence rate of the full discretization relies on the regularity of the driven noise process.
Particularly when the condition \eqref{eq:intro-AQ-condition} is fulfilled with $ \gamma = 2 $,  a classical convergence rate of order
$O(h^2 +\tau)$ for the backward Euler-finite element full discretization is reachable, even in multiple spatial dimensions
(see Remark \ref{rem:sectoin2-noise}). These findings are identical to those in \cite{wang2017strong},
where the strong convergence rate of the linear implicit Euler finite element scheme was analyzed for SPDEs with globally Lipschitz nonlinearity.
Once the nonlinearity grows super-linearly,  one can in general not expect the usual nonlinearity-explicit time-stepping schemes
that work well in the globally Lipschitz setting converge in the strong sense (see comments following Theorem 2 in \cite{jentzen2009pathwise}
and the relavant divergence result \cite{hutzenthaler2011strong}).
To address this issue, we therefore take the backward Euler, a nonlinearity-implicit scheme, for the temporal discretization.
Although some error estimates are taken from \cite{kruse2012optimal,kruse2014optimal,wang2017strong},
the presence of the non-globally Lipschitz (cubic) nonlinearity in the underlying model brings about essential difficulties in
the error analysis (see the proof of Theorems \ref{them-error-semi-discrete-problem}, \ref{them:error-estimate-full-FME})
and the error analysis  becomes much more involved than that in the globally Lipschitz SPDE setting.

In the following, we take error estimates of the spatial semi-discretization to illuminate our approach of the error analysis.
By introducing an auxiliary approximation process $\widetilde{X}_h$, defined by
\begin{align}\label{eq:intro-semi-auxiliary}
\widetilde{X}_h(t)
=
\mathcal{S}_h(t) P_h X_0
+
\int_0^t \mathcal{S}_h (t - s) P_hF(X(s))\,\dd s
+
\int_0^t \mathcal{S}_h ( t - s ) P_h \, \mathrm{d} W(s),
\quad
t \in [0, T],
\end{align}
we separate the spatial error  $\|X(t)-X_h(t)\|_{L^{2p}(\Omega;H)} $ into two parts,
\begin{equation}
\|X(t)-X_h(t)\|_{L^{2p}(\Omega;H)}
\leq
\| X(t) - \widetilde{X}_h(t) \|_{L^{2p}(\Omega;H)}
+
\| \widetilde{X}_h(t)  -  X_h(t) \|_{L^{2p}(\Omega;H)}.
\end{equation}
Subtracting \eqref{eq:intro-semi-auxiliary} from \eqref{eq:intro-mild-solution}, one can treat the first error term
directly and get $\| X(t) - \widetilde{X}_h(t) \|_{L^{2p}(\Omega;H)} = O(h^\gamma)$,
with the aid of existing estimates for the error operators  $\Psi_h(t):=\mathcal{S}(t)-\mathcal{S}_h(t)P_h$
and regularity properties of the mild solution $ X (t) $
(see estimates of $I_1, I_2, I_3, I_4$ in the proof of Theorem \ref{them-error-semi-discrete-problem} for details).
To bound the remaining error term $ \widetilde{e}(t) := \widetilde{X}_h(t)  -  X_h(t) $, we subtract
\eqref{eq:intro-FEM-mild-form} from \eqref{eq:intro-semi-auxiliary} to eliminate the stochastic convolution and  thus
$ \widetilde{e}(t) $ is time differentiable and satisfies
\begin{equation}
\tfrac{\dd}{\dd t} \widetilde{e}(t)+A_h\widetilde{e}(t)
=
P_h(F(X(t))-F(X_h(t))),
\quad
t \in (0 , T],
\quad
 \widetilde{e}_h(0)=0.
\end{equation}
Deterministic calculus together with the monotonicity of the nonlinearity, regularity properties of $\widetilde{X}_h(t)$,
the previous estimate of  $\| X(t) - \widetilde{X}_h(t) \|_{L^{2p}(\Omega;H)} $ and Gronwall's inequality facilitates the derivation of
$ \| \widetilde{e}(t) \|_{L^{2p}(\Omega;H)} = O (h^\gamma)$
(see \eqref{eq:estimate-2rd-error-term-semi}-\eqref{eq:error-discrete-and-auxiliary-control-by-auxiliary-error-semi}).

In the same manner as the semi-discrete case, we introduce an auxiliary  process
\begin{align}\label{eq:intro-full-auxiliary}
\begin{split}
\widetilde{X}_{h,m}
=
\mathcal{S}^m_{\tau,h}P_hX_0
+
\tau\sum_{i=0}^{m-1}\mathcal{S}^{m-i}_{\tau,h}P_hF(X(t_{i+1}))
+
\sum_{i=0}^{m-1} \mathcal{S}^{m-i}_{\tau,h}P_h \Delta W_{i+1},
\end{split}
\end{align}
and decompose the full discretization error  $ \| X(t_m) - X_{h,m} \|_{L^{2p}(\Omega;H)} $ as
\begin{equation}
\| X(t_m) - X_{h,m} \|_{L^{2p}(\Omega;H)}
\leq
\| X(t_m) - \widetilde{X}_{h,m} \|_{L^{2p}(\Omega;H)}
+
\| \widetilde{X}_{h,m}  -  X_{h,m} \|_{L^{2p}(\Omega;H)}.
\end{equation}
Then following the basic line as above, but with much more efforts made to exploit discrete versions of arguments
as used in the semi-discrete scenario, enables us to attain the desired error bounds for the full discretization
(cf. section \ref{sect:full-error}).


As usual, the strong convergence rate analysis of numerical SPDEs with super-linearly growing nonlinearities are
carried out based on appropriate uniform a priori moment $L_\infty$-bounds of approximations \cite{becker2017strong,liu2018strong,
feng2017finite,Majee2017optimal}.
Originally we develop a new approach of error analysis here, which does not rely on
high-order spatial regularity properties (e.g., a priori moment $L_\infty$-bounds) of approximation processes $X_h(t)$, $X_{h,m}$.
As already illustrated above, the new approach proposed for the error analysis is easy to understand
and can be extended to the error analysis for the stochastic Cahn-Hilliard equation \cite{kovacs2011finite,furihata2018strong,qi2018strong}.





Before closing the introduction part, we recall a few existing closely relevant works.
The backward Euler time semi-discretization was also examined in \cite{kovacs2015backward,kovacs2015discretisation}
for the problem \eqref{SACE}, with no spatial discretization.  Under assumption \eqref{eq:intro-AQ-condition} taking $ \gamma = 2 $,
i.e., $\|A^{\frac12}Q^{\frac12}\|_{\mathcal{L}_2}<\infty $,
only a strong convergence rate of order $\tfrac12$ was attained in \cite{kovacs2015discretisation}.
In \cite{brehier2018analysis,brehier2018strong,brehier2018weak}, pure time semi-discretizations of splitting type were studied
for \eqref{SACE}. Particularly, the authors of \cite{brehier2018strong} used exponential integrability properties of exact and
numerical solutions to identify a strong convergence rate of  order $1$, but only valid in one space dimension,
when the additive noise is moderately smooth, i.e., $ \gamma = 2 $ in \eqref{eq:intro-AQ-condition}.
Besides, various discretizations were investigated in \cite{becker2017strong,becker2016strong,liu2018strong,wang2018efficient}
 for the space-time white noise case
and in \cite{feng2017finite,Majee2017optimal} for the multiplicative one-dimensional noise case (gradient type noise in \cite{feng2017finite}),
only involved with a standard $\R$-valued Brownian motion.

The outline of this paper is as follows. In the next section, some preliminaries are collected and the
well-posedness and regularity properties of the considered problem are elaborated.
Section \ref{sect:Error-FEM} is devoted to error estimates of the finite element spatial semi-discretization
and section \ref{sect:full-error} provides error estimates of the backward Euler-finite element full discretization.
At the end of the article, some numerical examples are presented, illustrating the above theoretical findings.
\section{The stochastic Allen-Cahn equation}
%
%
Let $\mathbb{N} : = \{1,2,3,…\}$.
Given a separable $\mathbb{R}$-Hilbert space $\left(H,\left<\cdot,\cdot\right>,\|\cdot\|\right)$,
by $\mathcal{L}(H)$ we denote the Banach space of all linear bounded operators from $H$ into $H$.
Also, we denote by $ \mathcal{L}_2 (H)$  the Hilbert space consisting of Hilbert-Schmidt operators from $H$ into $H$, equipped with the {\color{black}{inner}}  product and the norm
\begin{align}
\left<\Gamma_1,\Gamma_2\right>_{\mathcal{L}_2(H)}
=
\sum_{i\in \mathbb{N}}\left<\Gamma_1\phi_i,\Gamma_2\phi_i\right>,
\;
\| \Gamma \|_{\mathcal{L}_2(H)}
=
\Big (
\sum_{i\in \mathbb{N}}\| \Gamma \phi_i\|^2
\Big)^\frac12
,
\end{align}
independent of the choice of orthonormal basis $\{\phi_i\}$ of $H$.
If $ \Gamma \in \mathcal{L}_2(H)$ and
$L\in \mathcal{L}(H)$, then $ \Gamma L$, $ L \Gamma \in \mathcal{L}_2(H)$ and
\begin{equation}
\| \Gamma L\|_{\mathcal{L}_2(H)}
\leq
\| \Gamma \|_{\mathcal{L}_2(H)}\|L\|_{\mathcal{L}(H)}
,
\;
 \|L \Gamma \|_{\mathcal{L}_2(H)}
 \leq
 \| \Gamma \|_{\mathcal{L}_2(H)}\|L\|_{\mathcal{L}(H)}.
\end{equation}
%
 %
 %
\subsection{Abstract framework and main assumptions}
In this subsection, we formulate main assumptions concerning the operators $A$ and  $Q$, the nonlinear term $F(\cdot)$, the noise term $W(t)$ and the initial value $X_0$, which will be used throughout this paper.

\begin{assumption}[Linear operator $A$]
\label{ass:equ-A-condition}
Let $\mathcal{D} \subset \R^d, d \in \{ 1, 2, 3 \} $ be a bounded open spatial domain with smooth boundary
and let $H : = L_2 ( \mathcal{D}, \R )$ be the real separable Hilbert space endowed with usual inner product
$\left<\cdot,\,\cdot\right>$ and the associated norm $\|\cdot\|=\left<\cdot,\,\cdot\right>^{\frac12}$.
Let $A:D(A)\subset H\rightarrow H$ be a densely defined, positive self-adjoint unbounded operator on $H$ with compact inverse.
\end{assumption}

Such assumptions imply the existence of a sequence of nondecreasing positive real numbers $\{\lambda_k\}_{k\geq1}$ and an orthonormal basis $\{e_k\}_{k\geq1}$ of $H$ such that
\begin{align} \label{eq:A-eigen-basis}
Ae_k=\lambda_ke_k,\quad \lim_{k\rightarrow\infty}\lambda_k=+\infty.
\end{align}
Furthermore,  it is known that $ - A$ generates an analytic semigroup  $\mathcal{S}(t) = e^{ - t A }$ satisfying
\begin{align}\label{spatial-temporal-S(t)}
\begin{split}
\|A^\mu \mathcal{S}(t)\|
&\leq
 Ct^{-\mu},\; t>0,\; \mu\geq 0,
 \\
 \|A^{-\nu}(I-\mathcal{S}(t))\|
 &\leq
 Ct^\nu,\quad t\geq0,\;\nu\in[0,1].
 \end{split}
 \end{align}
Throughout this article, we use generic constants which may vary at each appearance but are always independent of
discretization parameters.
By means of the spectral decomposition of $A$, we can also define the fractional powers $\gamma\in \mathbb{R}$ of $A$ in a simple way,
e.g., $A^\gamma v=\sum_{k=1}^\infty\lambda_k^\gamma \langle v,e_k \rangle e_k$. Then we denote the Hilbert space
$\dot{H}^\gamma : = D(A^{\frac \gamma2})$ with the inner product $\langle A^{\frac\gamma 2}\cdot, A^{\frac\gamma 2}\cdot\rangle$
and the associated norm $\| \cdot \|_{\gamma} : = \|A^{\frac \gamma 2}\cdot\|$.

\begin{assumption}[Nonlinearity]\label{assum:Nonlinearity}
Let $F:L_6(\mathcal{D};\mathbb{R})\rightarrow H$ be a deterministic mapping given by
\begin{align} \label{eq:F-f-Defn}
F(v)(x)=f(v(x))=v(x)-v^3(x),\;x\in \mathcal{D},\;v\in L_6(\mathcal{D};\mathbb{R}).
\end{align}
\end{assumption}
Here and below, by $L_{ r } ( \mathcal{D}; \R ), r \geq 1$ ($L_{ r } ( \mathcal{D} )$ or  $L_{ r }$ for short)
we denote a Banach space consisting of $r$-times integrable functions.
It is easy to check that, for any $v, \psi,  \psi_1, \psi_2 \in L_6(\mathcal{D};\mathbb{R}),$
\begin{equation}\label{eq:definition-derivative-F}
\begin{split}
\big ( F'(v) (\psi) \big) (x)
& =
f'(v(x)) \psi ( x )
=
(
1 - 3 v^2 ( x )
)
\psi ( x ),
\quad
x\in \mathcal{D},
 \\
 \big( F''(v) ( \psi_1, \psi_2 ) \big ) (x)
 & =
 f''(v(x))  \psi_1 ( x ) \psi_2 ( x )
 =
 -6v(x) \psi_1 ( x ) \psi_2 ( x ),
 \quad
 x\in \mathcal{D}.
\end{split}
\end{equation}
Moreover, {\color{black}{the following inequality holds}}
\begin{align}\label{eq:local-lipschitz-condition-F}
\big<u-v,F(u)-F(v)\big>
\leq
{\color{black}{\|u-v\|^2}},\;u,v\in L_6(\mathcal{D};\mathbb{R}).
\end{align}
\begin{assumption}[Noise process]
\label{assum:eq-noise}
Let $\{W(t)\}_{t\in[0,T]}$ be a standard $H$-valued $Q$-Wiener process on the stochastic basis
$\big(\Omega,\mathcal{F},\mathbb{P},\{\mathcal{F}_t\}_{t\in[0,T]}\big)$,
where the covariance operator $Q\in \mathcal{L}(H)$ is bounded, self-adjoint and positive semi-definite.
Assume
\begin{align}
\label{eq:A-Q-condition}
\|A^{\frac{\gamma-1}2 }Q^{\frac12}\|_{\mathcal{L}_2(H)}<\infty,
\quad
\text{ for some } \; {\color{black}{\gamma\in\big[\tfrac d 3,2\big]}},\;d\in\{1,2,3\}.
\end{align}
Additionally, we  assume that, for sufficiently large number $p_0\in \mathbb{N}$,
 \begin{align} \label{eq:WA-L18-condition}
\sup_{s\in[0,T]}\|W_A(s)\|_{L^{2p_0}(\Omega;L_{18})} < \infty,
\quad
\text{  with  }
\quad
W_{A}(t) : = \int_0^t\mathcal{S}(t-s)\,\mathrm{d} W(s).
 \end{align}
\end{assumption}

%
\begin{assumption}[Initial value]
\label{assum:initial-value-u0}
{\color{black}{For $\gamma \in [\tfrac{d}3,2]$ determined by the condition \eqref{eq:A-Q-condition}}}, we let  the initial data 
$X_0:\Omega\rightarrow H$ be $\mathcal{F}_0/\mathcal{B}(H)$-measurable and satisfy
\begin{align}
\mathbf{E}[\|X_0\|_\gamma^{p_0}]<\infty.
\end{align}
\end{assumption}
We remark that the assumption on the initial value can be relaxed, but at the expense of
having the constant $C$ later depending on $T^{-1}$, by {\color{black}{exploiting}} the
smoothing effect of the semigroup $ E (t), t \in [0, T]$ and standard non-smooth data error estimates \cite{thomee2006galerkin}.

{\color{black}{To conclude this subsection, we make some useful comments on Assumptions \ref{assum:eq-noise}.
\begin{remark}\label{rem:sectoin2-noise}
Note that the condition \eqref{eq:A-Q-condition} is commonly used in the literature \cite{kovacs2015discretisation,kovacs2015backward,kruse2014optimal,wang2017strong,yan2004semidiscrete}
but the condition \eqref{eq:WA-L18-condition} not.
Next we give more comments on when the unusual assumption \eqref{eq:WA-L18-condition} is fulfilled.
For the special case $Q=A^{-s}$, it is not difficult to see \eqref{eq:A-Q-condition} is fulfilled with $\gamma \in  [\frac d3, 2]$ iff
$ s > \gamma + \tfrac{d}{2} - 1 $, and  \eqref{eq:A-Q-condition} implies \eqref{eq:WA-L18-condition} since
$ \|W_A ( s ) \|_{L^{2p}(\Omega; L_{18})} \leq  \|W_A ( s ) \|_{L^{2p}(\Omega; C(\mathcal{D};\mathbb{R}) )} < \infty$
by \cite[Proposition 4.3]{da2012kolmogorov}. 
For the general case when $A$ and $Q$ do not own the same eigenbasis, we first recall the following 
Sobolev embedding inequalities, see, e.g.,
\cite[Theorem 7.57]{adams1975sobolev}  and  \cite[Lemma 3.1]{thomee2006galerkin},
\begin{align}\label{Sobolev-embedding-inequlaity-l18}
\dot{H}^\eta(\mathcal{D})\subset L_{18}(\mathcal{D}),
\quad
\mathcal{D} \subset \mathbb{R}^d, \,
d \in \{ 1, 2, 3\} 
\; \text{ for } \;\;\eta=\tfrac{4d}9.
\end{align}
This together with the later regularity estimate \eqref{lem:spatial-regularity-stochatic-convolution} implies 
that \eqref{eq:WA-L18-condition} can hold true for the general case,
provided \eqref{eq:A-Q-condition} is satisfied with $\gamma \geq \tfrac{4d}{9}$ in space dimension  $d\in\{1,2,3\}$. 
In particular, the space-time white noise ($Q = I$) in space dimension  $d = 1$ and the general trace-class noise ($\|Q^{\frac12}\|_{\mathcal{L}_2(H)}<\infty$) in space dimension  $d\in\{1,2\}$ are covered. In space dimension  $d=3$, smoother noise satisfying 
$ \|A^{\frac 1 6 }Q^{\frac12}\|_{\mathcal{L}_2(H)}<\infty $ is allowed.
Finally we emphasize that the assumption $ \gamma \geq \tfrac{d}{3}, d \in \{ 1, 2, 3 \} $ in \eqref{eq:A-Q-condition} is crucial in 
the following error estimates.
By contrast, the condition \eqref{eq:WA-L18-condition} is
only used to promise the well-posedness and the regularity estimate \eqref{eq:thm-wellposed-L6} of the nonlinear stochastic problem.
\end{remark}
}}
\subsection{Regularity results of the model}
In this part, we focus on the well-posedness of the underlying problem and the space-time regularity properties of the mild solution.
A preliminary theorem is stated as follows.
\begin{theorem}\label{thm:uniqueness-mild-solution}
Under Assumptions \ref{ass:equ-A-condition}-\ref{assum:initial-value-u0}, the problem \eqref{SACE} admits a unique mild solution, given by \eqref{eq:intro-mild-solution}, satisfying, for any $ p \geq 1 $,
\begin{equation}\label{eq:thm-wellposed-L6}
\sup_{s\in[0,T]}
\|X(s)\|_{L^{2p}(\Omega; L_6)}
\leq
C
\Big (
1 +
\|X_0 \|_{L^{2p}(\Omega; L_6)}
+
\sup_{ s \in [0, T] } \|W_A ( s ) \|_{L^{2p}(\Omega; L_{18})}^3
\Big ).
\end{equation}
\end{theorem}
To arrive at Theorem \ref{thm:uniqueness-mild-solution}, one can simply adapt the proof of \cite[Theorem 4.8]{da2012kolmogorov},
where the existence and uniqueness of the mild solution of the stochastic Allen-Cahn equation \eqref{SACE}
was established in the special case $Q=A^{-s}$.
There a basic tool for the proof is provided by the Yosida approximate arguments
and the assumption $Q=A^{-s}$, $s>\frac d2-1$ can be replaced by
assumptions {\color{black}{\eqref{eq:A-Q-condition}}} instead.
Indeed,  the assumption $Q=A^{-s}$, $s>\frac d2-1$ was simply used there to ensure
$ \|W_A ( s ) \|_{L^{2p}(\Omega; C(\mathcal{D};\mathbb{R}) )} < \infty$
and thus $ \|W_A ( s ) \|_{L^{2p}(\Omega; L_{18})} < \infty$.
Also, one can consult the proof of \cite[Theorem 5.5.8]{da1996ergodicity}.
%
{\color{black}{The}} above estimate \eqref{eq:thm-wellposed-L6} suffices to ensure
\begin{align}\label{eq:boundeness-f}
\sup_{s\in[0,T]}\|F(X(s))\|_{L^{2p}(\Omega;H)}
\leq
C
\big(
1
+
\sup_{s\in[0,T]}\|X(s)\|^{3}_{L^{6p}(\Omega;L_6)}
\big)
<\infty.
\end{align}
%
Equipped with Theorem \ref{thm:uniqueness-mild-solution}, we can get the following further regularity results.
\begin{theorem}\label{them:regulairty-mild-solution}
Under Assumptions \ref{ass:equ-A-condition}-\ref{assum:initial-value-u0},
the mild solution \eqref{eq:intro-mild-solution} enjoys the following regularity,
\begin{align}\label{them:spatial-regularity-mild-stoch}
\sup_{s\in[0,T]}
\|X(s)\|_{L^{2p} ( \Omega; \dot{H}^{\gamma} ) }
<
\infty
,
\quad
\forall p\geq 1,
\end{align}
and for any $\beta\in[0,\gamma]$,
\begin{align}\label{them:temporal-regularity-mild-stoch}
\| X(t) - X(s) \|_{L^{2p} ( \Omega; \dot{H}^{\beta} ) }
\leq
C(t-s)^{\frac{\min\{1,\gamma-\beta\}}2}.
\end{align}
\end{theorem}
Before proving Theorem \ref{them:regulairty-mild-solution}, we introduce some basic inequalities.
Recall first two well-known Sobolev embedding inequalities,
\begin{align}\label{eq:relation-C-H(deta)}
\dot{H}^{ \delta } \subset C(\mathcal{D};\mathbb{R}),
\;
\text{ for }
\;
\delta > \tfrac{d}{2},
\;
d \in
\{1,2,3\},
\end{align}
and
\begin{equation} \label{eq:H1-L6}
\dot{H}^1{\color{black}{\subset \dot{H}^{\frac d3}}}\subset L_6(\mathcal{D}),
\;
\text{ for }
\;
d \in
\{1,2,3\}.
\end{equation}
With \eqref{eq:relation-C-H(deta)} at hand, one can show
\begin{align}\label{eq:relation-L(1)-H(-deta))}
\begin{split}
\|A^{-\frac{\delta}2} x\|
&=
\sup_{\|v\|=1,v\in H}|\big<x, A^{-\frac {\delta} 2} v\big>|
\leq
\sup_{\|v\|=1,v\in H}\|x\|_{L_1}\|A^{-\frac{\delta} 2} v\|_{C(\mathcal{D};\mathbb{R})}
\\
&\leq
C\sup_{\|v\|=1,v\in H}\|x\|_{L_1}\|v\|
\leq
C
\|x\|_{L_1},
\quad
\forall \delta \in(\tfrac32, 2), x\in L_1(\mathcal{D}).
\end{split}
\end{align}
Similarly, but with the help of \eqref{eq:H1-L6}, one can find
\begin{align}\label{eq:relation-L(65)-H(-1)}
\begin{split}
\|A^{-\frac12}x\|
&
=
\sup_{\|\chi\|=1,\chi\in H}\big|\big<x, A^{-\frac12}\chi\big>\big|
\leq
\sup_{\|\chi\|=1,\chi\in H}
\|x\|_{L_{\frac65}} \|A^{-\frac12}\chi\|_{L_6}
\\
&\leq
C
\sup_{\|\chi\|=1,\chi\in H}
\|x\|_{L_{\frac65}} \|\chi\|
\leq
C
\|x\|_{L_{\frac65}},
\quad
\forall x\in L_{\frac65}(\mathcal{D}).
\end{split}
\end{align}
A slight modification of the proof of \cite[Theorem 3.1, Corollary 5.2]{kruse2012optimal}
gives the following lemma.
 \begin{lemma}\label{lem:regulairty-stochatic-convolution}
If condition \eqref{eq:A-Q-condition} from Assumption \ref{assum:eq-noise} is  valid,
then 
$\forall p\geq 1$,
\begin{align}\label{lem:spatial-regularity-stochatic-convolution}
\sup_{s\in[0,T]}
\|W_A(s)\|_{L^{2p} ( \Omega; \dot{H}^{\gamma} ) }
\leq
C\|A^{\frac{\mathbf{\gamma}-1}2}Q^{\frac12}\|_{\mathcal{L}_2(H)},
\end{align}
and, for any $\alpha\in[0,\gamma]$ and for $0\leq s<t\leq T$,
\begin{align}\label{lem:temporal-regularity-stochatic-convolution}
\|W_A(t)-W_A(s)\|_{L^{2p} ( \Omega; \dot{H}^{\alpha} ) }
\leq
C (t-s)^{\frac{\min\{1,\gamma-\alpha\}}2}
\| A^{\frac{\gamma-1}2}Q^{\frac12}\|_{\mathcal{L}_2(H)}.
\end{align}
\end{lemma}
In addition to the above preparations, we also need a lemma quoted from \cite[Lemma 3.2]{kruse2012optimal}.
\begin{lemma}
For any $\rho\in[0,1]$ and for all $x\in H$, it holds
 \begin{align}\label{spatial-temporal-integrand-S(t)}
 \begin{split}
 \int_{\tau_1}^{\tau_2}
 \|A^{\frac\rho2}\mathcal{S}(\tau_2-\sigma)x\|^2\,\dd \sigma
 &
 \leq
 C(\tau_2-\tau_1)^{1-\rho}\|x\|^2,\;0\leq\tau_1< \tau_2,
 \\
 \Big\|
 A^\rho\int_{\tau_1}^{\tau_2}
 \mathcal{S}(\tau_2-\sigma)x\,\dd \sigma
 \Big\|
 &
 \leq
 C(\tau_2-\tau_1)^{1-\rho}\|x\|,\;\;0\leq\tau_1< \tau_2.
 \end{split}
\end{align}
\end{lemma}
At the moment, we are able to start the proof of Theorem \ref{them:regulairty-mild-solution}.

{\it Proof of Theorem \ref{them:regulairty-mild-solution}.}
We take any fixed number $\delta_0\in(\frac32,2)$ and consider two possibilities:
either $\gamma \in {\color{black}{[\frac d3, \delta_0]}}$ or $\gamma \in (\delta_0, 2]$.
When \eqref{eq:A-Q-condition}  is fulfilled with $\gamma \in {\color{black}{[\frac d3, \delta_0]}}$,
we utilize \eqref{eq:boundeness-f}, \eqref{lem:spatial-regularity-stochatic-convolution} and \eqref{spatial-temporal-S(t)} with $\mu=\gamma$
to show, for $\gamma\in{\color{black}{[\frac d3,\delta_0]}}$,
\begin{align}\label{eq:bound-mild-solution-low-gamma}
\|X(t)\|_{L^{2p}(\Omega;\dot{H}^\gamma)}
\leq
&
\|
\mathcal{S}(t)X_0
\|_{L^{2p}(\Omega;\dot{H}^\gamma)}
+
\Big\|
\int_0^t\mathcal{S}(t-s) F(X(s))\,\dd s
\Big\|_{L^{2p}(\Omega;\dot{H}^\gamma)}
+
\|W_A(t)\|_{L^{2p}(\Omega;\dot{H}^\gamma)}
\nonumber \\
\leq&
C\|X_0\|_{L^{2p}(\Omega;\dot{H}^\gamma)}
+
C\int_0^t(t-s)^{-\frac\gamma2}\|F(X(s))\|_{L^{2p}(\Omega;H)}\,\dd s
+
C\|A^{\frac{\gamma-1}2}Q^{\frac12}\|_{\mathcal{L}_2}
\nonumber \\
\leq&
C\|X_0\|_{L^{2p}(\Omega;\dot{H}^\gamma)}
+
C
\sup_{s\in[0,T]}\|F(X(s))\|_{L^{2p}(\Omega;H)}
+
C\|A^{\frac{\gamma-1}2}Q^{\frac12}\|_{\mathcal{L}_2}
<\infty.
\end{align}
Concerning the temporal regularity of the mild solution, we apply \eqref{lem:spatial-regularity-stochatic-convolution}, \eqref{lem:temporal-regularity-stochatic-convolution}, \eqref{eq:boundeness-f} and \eqref{spatial-temporal-S(t)} 
 to obtain, for any $\beta\in[0,\gamma]$ {\color{black}{with $\gamma \in [\frac d3, \delta_0]$}},
\begin{align}
\|X(t)-X(s)\|_{L^{2p}(\Omega;\dot{H}^\beta)}
\leq &
\|(\mathcal{S}(t-s)-I)X(s)\|_{L^{2p}(\Omega;\dot{H}^\beta)}
+
\Big\|
\int_s^t\mathcal{S}(t-r)F(X(r))
\,\dd r
\Big\|_{L^{2p}(\Omega;\dot{H}^\beta)}
\nonumber \\
& +
\|W_A(t)-W_A(s)\|_{L^{2p}(\Omega;\dot{H}^\beta)}
+
\|(I-\mathcal{S}(t-s))W_A(s)\|_{L^{2p}(\Omega;\dot{H}^\beta)}
\nonumber \\
\leq &
C(t-s)^{\frac{\gamma-\beta}2}\sup_{s\in[0,T]}\|X(s)\|_{L^{2p}(\Omega;\dot{H}^\gamma)}
+
C\int_s^t(t-r)^{-\frac{\beta}2}\|F(X(r))\|_{L^{2p}(\Omega;H)}\,\dd r
\nonumber \\
&
+
C(t-s)^{\frac{\min\{1,\gamma-\beta\}}2}\|A^{\frac{\gamma-1}2}Q^{\frac12}\|_{\mathcal{L}_2}
+
C(t-s)^{\frac{\gamma-\beta}2}\|W_A(s)\|_{L^{2p}(\Omega;\dot{H}^\gamma)}
\nonumber \\
\leq &
C ( t - s )^{\frac{\min\{1,\gamma-\beta\}}2}
\Big(
\sup_{s\in[0,T]}\|X(s)\|_{L^{2p}(\Omega;\dot{H}^\gamma)}
\nonumber \\
&
\quad +
\sup_{s\in[0,T]}\|F(X(s))\|_{L^{2p}(\Omega;H)}
+
\|A^{\frac{\gamma-1}2}Q^{\frac12}\|_{\mathcal{L}_2}
\Big)
\nonumber \\
&
\leq
C(t-s)^{\frac{\min\{1,\gamma-\beta\}}2}.
\label{eq:temporal-regulairty-mild-solution-low-gamma}
\end{align}

Next, let us look at the other case $\gamma\in(\delta_0,2]$.
In this case, one can see
$\sup_{s\in[0,T]} \|X(s)\|_{L^{8p}(\Omega;\dot{H}^{\delta_0})}^{8p} < \infty$, as already verified in the former case.
Accordingly, applying \eqref{eq:relation-C-H(deta)} and \eqref{eq:temporal-regulairty-mild-solution-low-gamma} implies
\begin{align}
\|F(X(t))-F(X(r))\|_{L^{2p}(\Omega;H)}
&\leq
\big\|\,\|X(t)-X(r)\|(1+\|X(t)\|^2_{C(\mathcal{D},\mathbb{R})}+\|X(r)\|_{C(\mathcal{D},\mathbb{R})}^2)\big\|_{L^{2p}(\Omega;\mathbb{R})}
\nonumber \\
&\leq
\|X(t)-X(r)\|_{L^{4p}(\Omega;H)}
\big(
1 + \sup_{s\in[0,T]}\|X(s)\|^{2}_{L^{8p}(\Omega;\dot{H}^{\delta_0})}
\big)
\nonumber \\
&\leq
C|t-r|^{\frac12}.
\label{eq:temporall-regulairty-f(x)-high-gamma}
\end{align}
This together with \eqref{spatial-temporal-S(t)}, \eqref{spatial-temporal-integrand-S(t)} and \eqref{eq:boundeness-f} leads to, for $\beta\in[0,\gamma]$ {\color{black}{with $\gamma \in (\delta_0, 2]$}},
\begin{align}\label{eq:estiamte-semigroup-F-integrand}
\begin{split}
&
\Big\|
\int_s^t\mathcal{S}(t-r) F(X(r))\,\dd r
\Big\|_{L^{2p}(\Omega;\dot{H}^\beta)}
\\
& \quad \leq
\big\|\int_s^t\mathcal{S}(t-r)F(X(t))\,\dd r\big\|_{L^{2p}(\Omega;\dot{H}^\beta)}
+
\int_s^t\big\|\mathcal{S}(t-r)(F(X(t))-F(X(r)))\big\|_{L^{2p}(\Omega;\dot{H}^{\beta})}\,\dd r
\\
& \quad \leq
C(t-s)^{\frac{2-\beta}2}\|F(X(t))\big\|_{L^{2p}(\Omega;H)}
+
C
\int_s^t(t-r)^{-\frac\beta2}\|F(X(t))-F(X(r))\big\|_{L^{2p}(\Omega;H)}\,\dd r
\\
&
\quad
\leq
C(t-s)^{\frac{2-\beta}2}\sup_{s\in[0,T]}\|F(X(s))\big\|_{L^{2p}(\Omega;H)}
+
C\int_s^t(t-r)^{\frac{1-\beta}2}\,\dd r
\\
&
\quad
\leq
 C(t-s)^{\frac{2-\beta}2}.
\end{split}
\end{align}
Bearing this in mind and following the proof of \eqref{eq:bound-mild-solution-low-gamma}
and \eqref{eq:temporal-regulairty-mild-solution-low-gamma}, we can show \eqref{them:spatial-regularity-mild-stoch}
and \eqref{them:temporal-regularity-mild-stoch} in the case $\gamma\in(\delta_0,2]$.
The proof of Theorem \ref{them:regulairty-mild-solution} is thus complete.
$\square$
%


{\color{black}{
\begin{remark}
We would like to point out that the regularity estimate \eqref{eq:estiamte-semigroup-F-integrand} 
is only useful for the border case $\gamma = 2$, 
and the desired estimate for $\gamma< 2$ can be obtained directly by integration of the singularity $(t - r)^{-\frac\gamma 2}$.
Similar arguments are applied in the error analysis later (see \eqref{eq:estiamte-I3} and \eqref{eq:boundness-L22}),  
when handling the  border case $\gamma = 2$.
\end{remark}}}
\section{Error estimates of the spatial semi-discretization}
\label{sect:Error-FEM}
This section is devoted to error estimates of the finite element approximation of the stochastic problem \eqref{SACE}.
%
%
For the sake of simplicity, from here to section \ref{sect:full-error} we always assume that $\mathcal{D}\subset \mathbb{R}^d, d=1,2,3$, is an open convex polynomial domain and $A=-\Delta$ with $D(A)=H^2 (\mathcal{D}) \cap H_0^1(\mathcal{D})$.
%

In order to introduce the semi-discrete finite element approximation,
we present some notation and operators on the finite element space. Let $V_h\subset H_0^1(\mathcal{D})$, $h\in(0,1]$ be the space of
continuous functions that are piecewise linear over the triangulation $\mathcal{T}_h$  of $\mathcal{D}$.
Then we introduce a discrete Laplace operator $A_h:V_h\rightarrow V_h$ defined by
\begin{align}\label{eq:definition-discrete-A}
\left<A_hv_h,\chi_h\right>=a(v_h,\chi_h):=\big<\nabla v_h, \nabla \chi_h\big>,
\quad
\forall v_h,\;\chi_h\in V_h,
\end{align}
and a generalized projection operator $P_h:\dot{H}^{-1}\rightarrow V_h$ given by
\begin{align}\label{eq:defn-Ph}
\left<P_h v,\,\chi_h\right>=\left<v,\,\chi_h\right>,
\quad
\forall v\in \dot{H}^{-1},
\:
\chi_h\in V_h .
\end{align}
It is well-known that, the operators $A$ and $A_h$ obey
\begin{align}\label{eq:relation-A-Ah}
C_1\|A_h^{\frac r2}P_hv\|
\leq
\|A^{\frac r2}v\|
\leq
C_2
\|A_h^{\frac r2}P_hv\|,\;v\in \dot{H}^{r},\;r\in[-1,1].
\end{align}
The semi-discrete finite element method for the problem \eqref{SACE} is to find $X_h(t)\in V_h$ such that
\begin{align}\label{eq:semi-discrete-FEM}
\begin{split}
\dd X_h(t) +A_hX_h(t)\,\dd t=P_hF(X_h(t))\,\dd t+P_h\,\dd W(t),
\;
t \in (0, T],
\quad
X_h(0)=P_hX_0.
\end{split}
\end{align}
Let $\mathcal{S}_h(t)$ be the strongly continuous semigroup generated by the discrete Laplace operator $- A_h$.
Then it is easy to check that the semi-discrete problem \eqref{eq:semi-discrete-FEM} admits a unique solution in $V_h$,
given by $X_h(0)=P_hX_0$ and
\begin{align}
X_h(t)=\mathcal{S}_h(t)P_hX_0
+
\int_0^t\mathcal{S}_h(t-s)P_hF(X_h(s))\,\dd s
+
W_{A_h}(t),
\quad
t \in (0, T]
,
\end{align}
with
$W_{A_h}(t) : =\int_0^t\mathcal{S}_h(t-s)P_h\,\dd W(s)$.
The resulting spatial approximation error is measured  as follows.
\begin{theorem}\label{them-error-semi-discrete-problem}
Let $X(t)$ and $X_h(t)$ be the mild solutions of \eqref{SACE} and  \eqref{eq:semi-discrete-FEM}, respectively. If  Assumptions \ref{ass:equ-A-condition}-\ref{assum:initial-value-u0} are valid, then $ \forall p \in [1, \infty) $,
\begin{align}\label{them:error-results-seimidiscrete-problem}
\begin{split}
\|X(t)-X_h(t)\|_{L^{2p}(\Omega;H)} \leq
Ch^\gamma,
\quad
{\color{black}{\gamma \in[\tfrac d3,2]}}.
\end{split}
\end{align}
\end{theorem}
Its proof is postponed after we have been well-prepared with some important lemmas.
Define the semi-discrete approximation operator $\Psi_h(t), t \in [0, T]$ as follows,
\begin{equation}
\Psi_h(t) := \mathcal{S}(t) - \mathcal{S}_h(t)P_h, \quad t \in [0, T].
\end{equation}
The following results listed in \cite[Lemmas 4.1, 4.2]{kruse2014optimal} on the error operator $\Psi_h(t)$ are crucial
in the error estimates of the semi-discrete finite element approximation.
\begin{lemma}\label{lem:error-estimates-semi-determi}
Under  Assumption \ref{ass:equ-A-condition}, the following estimates for the error operator
$\Psi_h(t)$ hold.\\
(i) For $0\leq \nu\leq \mu\leq 2$, it holds that
\begin{align}\label{lem:error-determin-initial-spatial-step-positive}
\|\Psi_h(t)x\|
\leq Ch^\mu t^{-\frac{\mu-\nu}2}\|x\|_\nu,\;for\;all\;x\in \dot{H}^\nu, \;t>0.
\end{align}
(ii) Let $0\leq\rho\leq1$. Then
\begin{align}\label{lem:error-derterminsitic-integrand}
\Big\|
\int_0^t\Psi_h(s)x\,\mathrm{d} s
\Big\|
\leq
Ch^{2-\rho}\|x\|_{-\rho},\;
for\;all\;x\in \dot{H}^{-\rho},\;t>0.
\end{align}
(iii) Let $0\leq\varrho\leq1$. Then
\begin{align}\label{lem-error-determin-integrand-square-norm}
\Big(\int_0^t\left\|\Psi_h(s)x\right\|^2\,\mathrm{d} s\Big)^{\frac12}
\leq
Ch^{1+\varrho}\|x\|_{\varrho},\;
for\;all\;x\in \dot{H}^{\varrho},\;t>0.
\end{align}
\end{lemma}
Additionally, we need smoothing properties of the semigroup $\mathcal{S}_h(t)$, as described below.
\begin{lemma}\label{lem:boundess-discrete-semigroup-integrand}
Under  Assumption \ref{ass:equ-A-condition}, the following estimates  for the discrete semigroup $\mathcal{S}_h(t)$ hold,
\begin{align}
\|A_h^{\frac\mu2}\mathcal{S}_h(t)P_hx\|
&\leq
Ct^{-\frac\mu2}\|x\|,
\quad
\forall
\mu\in[0,1],
\:
x\in H,
\label{lem:eq-boundess-discrete-semigroup}
\\
\int_0^t\|A_h^{\frac12}\mathcal{S}_h(s)P_hx\|^2\,\mathrm{d} s
&\leq
 C\|x\|^2,
\quad
\forall x \in H.
 \label{lem:eq-boundess-discrete-semigroup-integrand}
\end{align}
\end{lemma}
%
%
The assertion \eqref{lem:eq-boundess-discrete-semigroup} is obvious and \eqref{lem:eq-boundess-discrete-semigroup-integrand}
is derived by using (4.19) in \cite{yan2004semidiscrete}.
\begin{lemma}\label{lem:eq-spatial-semi-discrete-problem}
Suppose  Assumptions \ref{ass:equ-A-condition}-\ref{assum:initial-value-u0} hold. {\color{black}{Let $X_h(t)$ be the solution of \eqref{eq:semi-discrete-FEM} and denote
$Y_h(t)  :=X_h(t)-W_{A_h}(t)$, with
$W_{A_h}(t) : =\int_0^t\mathcal{S}_h(t-s)P_h\,\dd W(s)$.  Then
\begin{align}\label{lem:eq-bound-semi-problem}
\sup_{s\in[0,T]}\mathbf{E}\big[\|Y_h(s)\|^{2p}\big]
+
\int_0^T\mathbf{E}\big[\|\nabla Y_h(s)\|^2\big]\,\mathrm{d} s
<\infty.
\end{align}}}
\end{lemma}
{\it Proof of Lemma \ref{lem:eq-spatial-semi-discrete-problem}.}
Recall first that
$
Y_h(t)  = \mathcal{S}_h(t)P_hX_0
+
\int_0^t\mathcal{S}_h(t-s)P_hF(X_h(s))\,\dd s.
$
Then $Y_h(t)$ is time differentiable and obeys
\begin{align}\label{ASCE-semi-scheme-equival}
\tfrac{\dd}{\dd t}Y_h(t)+A_hY_h(t)
=
P_hF(Y_h(t)+W_{A_h}(t)),
\quad
Y_h(0)=P_hX_0.
\end{align}
By multiplying  both sides of \eqref{ASCE-semi-scheme-equival} by $Y_h(t)$, taking the  inner product and using \eqref{eq:local-lipschitz-condition-F}, we obtain
\begin{align}\label{regulairty-mild-solution-semi-SACE-equvien-in-L2-norm}
\tfrac12\tfrac{\dd }{\dd s}\|Y_h(s)\|^2  +  \|\nabla Y_h(s)\|^2
&=
\left<F(Y_h(s)+W_{A_h}(s))-F(W_{A_h}(s)),Y_h(s)\right> + \left<F(W_{A_h}(s)),Y_h(s)\right>
\nonumber \\
&\leq
C
\|Y_h(s)\|^2 + \tfrac12\|F(W_{A_h}(s))\|^2  +  \tfrac12\|Y_h(s)\|^2
\nonumber \\
&\leq
C \|F(W_{A_h}(s))\|^2  +  C \|Y_h(s)\|^2,
\end{align}
which, after integration over $[0, t]$ and using the Gronwall inequality, gives that
\begin{align}\label{eq:bound-semi-random-equaiton}
\|Y_h(t)\|^2+\int_0^t\|\nabla Y_h(s)\|^2\,\dd s
\leq
C
\Big(
\int_0^t\|F(W_{A_h}(s))\|^2\,\dd  s+\|P_hX_0\|^2
\Big).
\end{align}
Then, using \eqref{eq:F-f-Defn},  \eqref{eq:A-Q-condition}, \eqref{eq:H1-L6}, \eqref{eq:relation-A-Ah}, \eqref{lem:eq-boundess-discrete-semigroup-integrand} and the Burkholder-Davis-Gundy-type inequality shows
{\color{black}{\begin{align}\label{eq:boundness-f(WAH)}
\begin{split}
\sup_{s\in[0,T]} \|F(W_{A_h}(s))\|_{L^{2p}(\Omega;H)}
&
\leq
C(
1
+
\sup_{s\in[0,T]} \| W_{A_h} ( s ) \|_{L^{ 6 p } ( \Omega; L_6 ) }^{3}
)
\\
&
\leq
C(
1
+
\sup_{s\in[0,T]} \| W_{A_h} ( s ) \|_{L^{ 6 p } ( \Omega; \dot{H}^{d/3} ) }^{3}
)
\\
& \leq
C
\Big(
1+\sup_{s\in[0,T]}
\Big\|\int_0^sA_h^{\frac d 6}\mathcal{S}_h(s-r)P_h\,\dd W(r)\Big\|_{L^{ 6 p } ( \Omega; H ) }^{3}
\Big)
\\
&
\leq
C
\Big(
1 + \sup_{s\in[0,T]}
\Big(
\int_0^s
\big \|
A_h^{\frac d6}\mathcal{S}_h(s-r)P_hQ^{\frac12}
\big\|_{\mathcal{L}_2}^2\,\dd r
\Big)^{3/2}
\Big)
\\
&
\leq
C(1+\|A_h^{\frac{d-3}6}P_hQ^{\frac12}\|^{3}_{\mathcal{L}_2})
\\
&\leq
{\color{black}{C(1+\|A^{\frac{d-3}6}Q^{\frac12}\|^{3}_{\mathcal{L}_2})}}
<
\infty.
\end{split}
\end{align}}}
This combined with Assumption \ref{assum:initial-value-u0}  shows the desired assersion.
$\square$
We are now ready to {\color{black}{ prove}} Theorem \ref{them-error-semi-discrete-problem}.

{\it Proof of Theorem \ref{them-error-semi-discrete-problem}.}
By introducing the following auxiliary  process,
\begin{align}\label{eq:definition-mild-solution-semi-auxiliary}
\widetilde{X}_h(t)
=
\mathcal{S}_h(t)P_hX_0
+
\int_0^t\mathcal{S}_h(t-s)P_hF(X(s))\,\dd s
+
W_{A_h}(t),
\end{align}
we separate the considered error term $\|X(t)-X_h(t)\|_{L^{2p}(\Omega;H)}$ as
\begin{equation}\label{eq:spatial-error-decomp}
\|X(t)-X_h(t)\|_{L^{2p}(\Omega;H)}
\leq
\| X(t) - \widetilde{X}_h(t) \|_{L^{2p}(\Omega;H)}
+
\| \widetilde{X}_h(t)  -  X_h(t) \|_{L^{2p}(\Omega;H)}.
\end{equation}
In view of \eqref{eq:boundeness-f}, \eqref{eq:H1-L6},  \eqref{eq:relation-A-Ah} and \eqref{lem:eq-boundess-discrete-semigroup},
we acquire that, for any $ t \in [ 0, T] $,
{\color{black}{\begin{align}\label{eq:boundness-semi-discrete-auxiliary-problem}
\begin{split}
&\|\widetilde{X}_h( t )-W_{A_h}(t)\|_{L^{2p}(\Omega;L_6)}
\leq
\|\widetilde{X}_h( t )-W_{A_h}(t)\|_{L^{2p}(\Omega;\dot{H}^{\frac d3})}
\\
&\leq
\|\mathcal{S}_h(t)P_hX_0\|_{L^{2p}(\Omega;\dot{H}^{\frac d3})}
+
\int_0^t\|\mathcal{S}_h(t-s)P_hF(X(s))\|_{L^{2p}(\Omega;\dot{H}^{\frac d3})}\,\dd s
\\
&\leq
C\|X_0\|_{L^{2p}(\Omega;\dot{H}^{\frac d3})}
+
C
\sup_{s\in[0,T]}
\|F(X(s))\|_{L^{2p}(\Omega;H)}
\int_0^t{\color{black}{(t-s)^{-\frac d6}}}\,\dd s
<\infty.
\end{split}
\end{align}
Noting that $\| W_{A_h}(t) \|_{L^{2p}(\Omega;L_6)} < \infty$,  as implied by \eqref{eq:boundness-f(WAH)},
we know that
\begin{align}
\|\widetilde{X}_h( t )\|_{L^{2p}(\Omega;L_6)}<\infty.
\end{align}
}}

With this we start to bound the {\color{black}{first error term}} in \eqref{eq:spatial-error-decomp}.
Subtracting \eqref{eq:definition-mild-solution-semi-auxiliary} from \eqref{eq:intro-mild-solution} yields
\begin{align}
\|X(t)-\widetilde{X}_h(t)\|_{L^{2p}(\Omega;H)}
\leq
&
\|(\mathcal{S}(t)-\mathcal{S}_h(t)P_h)X_0\|_{L^{2p}(\Omega;H)}
\nonumber
\\
& +
\big\|\int_0^t(\mathcal{S}(t-s)-\mathcal{S}_h(t-s)P_h)F(X(t))\,\dd s\big\|_{L^{2p}(\Omega;H)}
\nonumber
\\
& +
\int_0^t\|(\mathcal{S}(t-s)-\mathcal{S}_h(t-s)P_h)(F(X(t))-F(X(s)))\|_{L^{2p}(\Omega;H)}\,\dd s
\nonumber
\\
& +
\big\|\int_0^t(\mathcal{S}(t-s)-\mathcal{S}_h(t-s)P_h)\,\dd W(s)\big\|_{L^{2p}(\Omega;H)}
\nonumber
\\
: =
&
I_1+I_2+I_3+I_4.
\end{align}
Subsequently $I_1, I_2, I_3$ and $I_4$ will be treated separately.  For the first term $I_1$, we utilize \eqref{lem:error-determin-initial-spatial-step-positive} with $\mu=\nu=\gamma$ to derive
\begin{align}\label{eq:estiamte-I1}
\begin{split}
I_1\leq Ch^\gamma\|X_0\|_{L^{2p}(\Omega;\dot{H}^\gamma)}.
\end{split}
\end{align}
Employing \eqref{eq:boundeness-f} and \eqref{lem:error-derterminsitic-integrand} with $\rho=0$ enables us to obtain
\begin{align}
I_2\leq
Ch^2\|F(X(t))\|_{L^{2p}(\Omega;H)}
\leq
Ch^2.
\end{align}
To handle $I_3$, we recall \eqref{eq:temporall-regulairty-f(x)-high-gamma} and \eqref{eq:boundeness-f}, which together imply,
for any fixed number $\delta_0\in(\frac32,2)$,
\begin{align}\label{prop:boundness-f(x)-f(y)-H}
\|F(X(t) ) - F(X(s))\|_{L^{2p}(\Omega;H)}
\leq
\left\{\begin{array}{ll}
C,& {\color{black}{\gamma\in[\frac d3,\delta_0]}},
\\
C | t - s |^{ \frac12 },& \gamma\in(\delta_0,2].
\end{array}\right.
\end{align}
Therefore,  using \eqref{lem:error-determin-initial-spatial-step-positive} with $\mu=\gamma$, $\nu=0$
and also taking \eqref{prop:boundness-f(x)-f(y)-H} into consideration result in
\begin{align}\label{eq:estiamte-I3}
I_3
\leq
Ch^\gamma\int_0^t (t-s)^{-\frac\gamma2}\|F(X(t))-F(X(s))\|_{L^{2p}(\Omega;H)}\,\dd s
\leq Ch^\gamma.
\end{align}
Now it remains to bound $I_4$. Combining the Burkholder-Davis-Gundy type inequality  and
\eqref{lem-error-determin-integrand-square-norm} with $\varrho=\gamma-1$ results in
\begin{align}
I_4
\leq
C_{p}
\Big(
\int_0^t\|(\mathcal{S}(t-s)-\mathcal{S}_h(t-s)P_h)Q^{\frac12}\|^2_{\mathcal{L}_2}\,\dd s
\Big)^{\frac12}
\leq
Ch^\gamma\|A^{\frac{\gamma-1}2}Q^{\frac12}\|_{\mathcal{L}_2}.
\end{align}
Finally, putting the above estimates together gives
\begin{align}\label{eq:error-X-widetilde(X)}
\|X(t)-\widetilde{X}_h(t)\|_{L^{2p}(\Omega;H)}
\leq
Ch^\gamma.
\end{align}
Next we turn our attention  to the error $\widetilde{e}(t):=\widetilde{X}_h(t)-X_h(t)$, which is time differentiable and
\begin{align}\label{eq:pde-e(t)}
\tfrac{\dd}{\dd t} \widetilde{e}(t)+A_h\widetilde{e}(t)
=
P_h(F(X(t))-F(X_h(t))),
\quad
t \in (0, T],
\quad
\widetilde{e}_h(0)=0.
\end{align}
Note that Lemma \ref{lem:eq-spatial-semi-discrete-problem} and \eqref{eq:boundness-semi-discrete-auxiliary-problem}
guarantee that $\sup_{s \in [0, T]} \mathbf{E} [\|\widetilde{e}(s)\|^{2p} ]  +  \int_0^T\mathbf{E} [\|\nabla \widetilde{e}(s) \|^2 ] \dd s < \infty$.
Multiplying both sides of \eqref{eq:pde-e(t)} by $\widetilde{e} (t)$,
and applying \eqref{eq:local-lipschitz-condition-F}, \eqref{eq:relation-L(65)-H(-1)}, \eqref{eq:defn-Ph} and \eqref{eq:definition-discrete-A}
 tell us
\begin{align} \label{eq:estimate-2rd-error-term-semi}
\begin{split}
&
\tfrac12\tfrac{\dd }{\dd s}\|\widetilde{e}(s)\|^2 + \left<\nabla \widetilde{e}(s),\nabla \widetilde{e}(s)\right>
\\
&
\quad =
\langle F(\widetilde{X}_h(s))-F(X_h(s)),\widetilde{e}(s)\rangle
+
\langle F(X(s))-F(\widetilde{X}_h(s)),\widetilde{e}(s) \rangle
\\
&
\quad \leq
C \|\widetilde{e}(s)\|^2
+
\|A^{-\frac12} (F(X(s))-F(\widetilde{X}_h(s)))\|\,\|\nabla \widetilde{e}(s)\|
\\
&
\quad \leq
C \|\widetilde{e}(s)\|^2
+
\tfrac12\|F(X(s))-F(\widetilde{X}_h(s))\|_{L_{\frac65}}^2
+
\tfrac12\|\nabla \widetilde{e}(s)\|^2.
\end{split}
\end{align}
Then  integrating over $[0,t]$ and using H\"{o}lder's inequality {\color{black}{give}} that
\begin{align}
\begin{split}
\|\widetilde{e}(t)\|^2
&\leq
C\int_0^t\|\widetilde{e}(s)\|^2\,\dd s
+
C
\int_0^t\|F(X(s))-F(\widetilde{X}_h(s))\|_{L_{\frac65}}^2\,\dd s
\\
&
\leq
C\int_0^t\|\widetilde{e}(s)\|^2\,\dd s
+
C\int_0^t\|X(s)-\widetilde{X}_h(s)\|^2\big(1+\|X(s)\|_{L_6}^4 + \|\widetilde{X}_h(s)\|_{L_6}^4\big)\,\dd s.
\\
\end{split}
\end{align}
Using Gronwall's inequality before employing \eqref{eq:boundness-semi-discrete-auxiliary-problem}, \eqref{eq:error-X-widetilde(X)}
and Theorem \ref{thm:uniqueness-mild-solution}, one can arrive at
\begin{align}\label{eq:error-discrete-and-auxiliary-control-by-auxiliary-error-semi}
\begin{split}
\| \widetilde{e}(t) \|_{L^{2p} (\Omega; H)}^2
&
\leq
C
\int_0^t
\|
X(s)-\widetilde{X}_h(s)
\|_{L^{4p} (\Omega; H)}^2
\big(
1
+
\|X(s)\|_{L^{8p} (\Omega; L_6)}^4
+
 \|\widetilde{X}_h(s)\|_{L^{8p} (\Omega; L_6)}^4
\big)
 \,\dd s
\\
&
\leq
Ch^{2\gamma},
\end{split}
\end{align}
which in a combination with \eqref{eq:error-X-widetilde(X)} shows \eqref{them:error-results-seimidiscrete-problem}, as required.
$\square$
%
%
\section{Error estimates of the spatio-temporal full discretization}
\label{sect:full-error}
In the present section, we proceed to study a full discretization based on the finite element semi-discretization.
Let $\tau:=T/M$, $M\in \mathbb{N}$ be a uniform time-step size and write $t_m=m\tau$, for $m\in\{1,2\cdots, M\}$.
We discrete \eqref{eq:semi-discrete-FEM} in time with a backward Euler scheme and the resulting {\color{black}{fully}} discrete problem is to 
find $\mathcal{F}_{t_m}$-adapted $V_h$-valued random variables $X_{h,m}, m \in \{1,2\cdots, M\}$ such that, 
\begin{align}\label{eq:full-discrete-method0}
\begin{split}
X_{h,m} =  X_{h,m-1}  - \tau A_h  X_{h,m}   + \tau P_hF(X_{h,m}) + P_h \Delta W_m,
\quad
X_{h,0}=P_hX_0,
\quad
m\in\{1,2,\cdots, M\},
\end{split}
\end{align}
or equivalently,
\begin{align}\label{eq:full-discrete-method}
\begin{split}
X_{h,m} = \mathcal{S}_{\tau,h} X_{h,m-1} + \tau \mathcal{S}_{\tau,h} P_hF(X_{h,m}) + \mathcal{S}_{\tau,h} P_h \Delta W_m,
\quad
X_{h,0}=P_hX_0,
\quad
m\in\{1,2,\cdots, M\},
\end{split}
\end{align}
where we write $\Delta W_m := W(t_m)-W(t_{m-1})$, $\mathcal{S}_{\tau,h}:=( I + \tau A_h)^{-1}$ for brevity.
{\color{black}{Observe that the time-stepping scheme \eqref{eq:full-discrete-method} is implicit in the nonlinear term.
The first main issue concerns the well-posedness of the scheme, which is addressed by 
Proposition \ref{prop:well-posedness-backward-Euler} below. 
To implement the time-stepping scheme in the numerical experiment later, we simply used the fixed point iteration to
obtain approximation solutions to the nonlinear implicit systems.}}
{\color{black}{
\begin{prop}[Well-posedness of the fully discrete scheme]
\label{prop:well-posedness-backward-Euler}
Let Assumptions \ref{ass:equ-A-condition}-\ref{assum:initial-value-u0} hold and let $ \tau \leq 1 $. 
The fully discrete scheme \eqref{eq:full-discrete-method0} (or  \eqref{eq:full-discrete-method}) 
has a unique solution $ \{ X_{h,m} \} _{ m \in \{1,2\cdots, M\} }$ 
in $V_h$, which is $\mathcal{F}_{t_m}$-adapted. 
\end{prop}
{\it Proof of Proposition \ref{prop:well-posedness-backward-Euler}.}
For $h > 0$ and $\tau \leq 1$ fixed, we define a function $G_{h, \tau} \colon V_h \rightarrow V_h$ on the finite dimensional space $V_h$, by
$G_{h, \tau} ( z) : = z + \tau A_h z - \tau P_h F (z), \, z \in V_h$. 
In the light of \eqref{eq:A-eigen-basis}, \eqref{eq:local-lipschitz-condition-F}, \eqref{eq:definition-discrete-A}, 
properties of $P_h$ and the assumption $ \tau \leq 1 $,
it is not difficult to check that $G_{h, \tau}$ is continuous in $V_h$ and
\begin{equation}
\langle G_{h, \tau} (z_1) - G_{h, \tau} (z_2), z_1 - z_2 \rangle \geq ( 1 + \lambda_{1} \tau - \tau ) \| z_1 - z_2 \|^2
\geq  \lambda_{1}  \tau  \| z_1 - z_2 \|^2, 
\quad
z_1, z_2 \in V_h,
\end{equation}
where we used $\lambda_{1} $ to mean the first eigenvalue of $A$.
%
Thanks to \cite[Theorem C.2]{stuart1996dynamical},  the implicit equation $ G_{h, \tau} ( z ) = b $ for any $b \in V_h$ admits a unique solution $z = G_{h, \tau}^{ -1 } ( b )$ in $ V_h $. This implies the well-posedness of the fully discrete scheme \eqref{eq:full-discrete-method0}, as required. 
$\square$
}
}  

Further,  the recurrence \eqref{eq:full-discrete-method} promises
\begin{align}\label{eq:solution-full-discete-method}
\begin{split}
X_{h,m}
&=
\mathcal{S}^m_{\tau,h}X_{h,0}
+
\tau\sum_{i=0}^{m-1}\mathcal{S}^{m-i}_{\tau,h}P_hF(X_{h,i+1})
+
W_{A_h}^m,
\quad
\text{ with }
\:
W_{A_h}^m :=
\sum_{i=0}^{m-1} \mathcal{S}^{m-i}_{\tau,h}P_h \Delta W_{i+1}.
\end{split}
\end{align}

\begin{theorem}\label{them:error-estimate-full-FME}
Let $X(t)$  be the mild solution of \eqref{SACE} and let $X_{h,m}$ be produced by \eqref{eq:full-discrete-method}.
If  Assumptions \ref{ass:equ-A-condition}-\ref{assum:initial-value-u0} are valid and {\color{black}{$\tau \leq \tfrac13$}}, then it holds that
\begin{align}\label{them:error-result-full-FME}
\|X(t_m)-X_{h,m}\|_{L^{2p}(\Omega;H)}
\leq
  C(h^\gamma+ \tau^{\frac\gamma2}),
  \quad
  {\color{black}{\gamma\in[\tfrac d 3,2].}}
\end{align}
\end{theorem}
Its proof is also postponed.
Define the fully discrete approximation operators $\Psi_{\tau,h}(t), t \in [0, T]$ as
\begin{equation}
\Psi_{\tau,h}(t)=\mathcal{S}(t)-\mathcal{S}_{\tau, h}^mP_h,
\quad
\forall \, t\in [t_{m-1},t_m), \,
m \in \{ 1,2,..., M \}.
\end{equation}
The forthcoming two lemmas, coming from \cite[Lemmas 4.3, 4.4]{kruse2014optimal},
are a temporal version of Lemmas \ref{lem:error-estimates-semi-determi},\ref{lem:boundess-discrete-semigroup-integrand},
and play a significant role in the error estimates of the full-discrete approximation.
\begin{lemma}
Under  Assumption \ref{ass:equ-A-condition}, the following estimates hold.\\
(i) For $0\leq\nu\leq\mu\leq 2$, it holds that
\begin{align}\label{eq:error-estimates-deter-full-intial-positive-norm}
\|\Psi_{\tau,h}(t)x\|
\leq
C(h^\mu+\tau^{\frac\mu2})t^{-\frac{\mu-\nu}2}\|x\|_{\nu}, \quad \text { for all } \,  x\in \dot{H}^\nu.
\end{align}
(ii) For $0\leq\rho\leq 1$,  it holds that
\begin{align}\label{lem:eq-error-full-sum}
\Big\|
\int_0^t \Psi_{\tau,h} (s)x\,\mathrm{d} s
\Big\|
\leq
C(h^{2-\rho}+\tau^{\frac{2-\rho}2}) \|x\|_{-\rho},\quad \text { for all } \, x\in \dot{H}^{-\rho}.
\end{align}
(iii) For $0\leq\varrho\leq 1$, it holds that
\begin{align}\label{lem:eq-error-determ-full-sum-square-norm}
\Big(
\int_0^t\| \Psi_{\tau,h} ( s )x\|^2\,\mathrm{d} s
\Big)^{\frac12}
\leq
C(h^{1+\varrho}+\tau^{\frac{1+\varrho}2}) \|x\|_{\varrho}, \quad \text { for all } \,  x\in \dot{H}^{\varrho}.
\end{align}
\end{lemma}
%
\begin{lemma}\label{lem:boundness-full-discrete}
Under Assumption \ref{ass:equ-A-condition}, the following estimates for $\mathcal{S}_{\tau,h}^m$ hold, for any $x\in H$
\begin{align}
\|A_h^{\frac\mu2}\mathcal{S}_{\tau,h}^mP_hx\|
&\leq
Ct_m^{-\frac\mu2}\|x\| ,
\;\mu\in[0,1],
\label{eq:boundness-full-discrete-operator}
\\
\tau\sum_{i=1}^m\|A_h^{\frac12}\mathcal{S}_{\tau,h}^iP_hx\|^2
&\leq
C\|x\|^2.
\label{eq:boundness-full-operator-sum}
\end{align}
\end{lemma}
%
%
\begin{lemma}\label{lem:boundess-full-problem}
Suppose  Assumptions \ref{ass:equ-A-condition}-\ref{assum:initial-value-u0} hold and {\color{black}{$\tau \leq \tfrac13$. 
Let $X_{h,m}$ be produced by \eqref{eq:full-discrete-method} and denote
$Y_{h,m}  := X_{h,m} - W_{A_h}^m$ with $W_{A_h}^m$ defined as in \eqref{eq:solution-full-discete-method}.}} 
Then 
\begin{align}\label{eq:boundness-full-FEM}
\sup_{ M \in \N }
\sup_{m\in\{1,2,\cdots, M\}}
\Big (
\mathbf{E}\left[\|Y_{h,m}\|^{2p}\right]
+
\tau\sum_{i=1}^m \mathbf{E}\big[\|\nabla Y_{h,i}\|^2\big]
\Big)
<\infty.
\end{align}
\end{lemma}
{\it Proof of Lemma \ref{lem:boundess-full-problem}.}
{\color{black}{Note first that $Y_{h,m}$ satisfies}}
\begin{equation}
Y_{h,m}
=
\mathcal{S}^m_{\tau,h}Y_{h,0}
+
\tau\sum_{i=0}^{m-1}\mathcal{S}^{m-i}_{\tau,h}P_hF( Y_{h,i+1} + W_{A_h}^{i + 1} ),
\;
Y_{h,0}=P_hX_0,\;
m\in\{1,2,\cdots, M\}
.
\end{equation}
It is {\color{black}{straightforward}}  to verify that $Y_{h,m}$ satisfies
\begin{align}
\tfrac{Y_{h,m}-Y_{h,m-1}}{\tau}+A_hY_{h,m}=P_hF(Y_{h,m}+W_{A_h}^m),
\quad
Y_{h,0}=P_hX_0,
\;
m\in\{1,2,\cdots, M\}.
\end{align}
Multiplying this equation by $Y_{h,m}$ and using \eqref{eq:local-lipschitz-condition-F}, \eqref{eq:definition-discrete-A} imply
\begin{align}
\begin{split}
& \left<Y_{h,m}-Y_{h,m-1},Y_{h,m}\right>
+
\tau\left<\nabla Y_{h,m},\nabla Y_{h,m}\right>
\\
&
\quad
=
\tau \big \langle F(Y_{h,m}+W_{A_h}^m)-F(W_{A_h}^m), Y_{h,m} \big \rangle
+
\tau \big \langle F(W_{A_h}^m), Y_{h,m}  \big \rangle
\\
&
\quad
\leq
\color{black}{\tfrac {9\tau}{8}}
\|Y_{h, m}\|^2
+
2 \tau \|F(W_{A_h}^m)\|^2.
\end{split}
\end{align}
Further, using the fact
$\frac12(\|Y_{h,m}\|^2-\|Y_{h,m-1}\|^2)
\leq
\left<Y_{h,m}-Y_{h,m-1},Y_{h,m}\right>$
 and summation on $m$ shows
 \begin{align}
\tfrac12\|Y_{h,m}\|^2
+
\tau\sum_{i=1}^m\|\nabla Y_{h,i}\|^2
\leq
\tfrac12\|Y_{h,0}\|^2
+
\tfrac {9\tau}{8} \sum_{i=1}^m\|Y_{h,i}\|^2
+
2 \tau\sum_{i=1}^m\|F(W_{A_h}^i)\|^2,
 \end{align}
{\color{black}{
which, after rearrangement and noting $ \tau \leq \tfrac13$, shows
\begin{equation}
\begin{split}
\tfrac14 \|Y_{h,m}\|^2 
+
& 
2 \tau\sum_{i=1}^m\|\nabla Y_{h,i}\|^2
 \leq
( 1 - \tfrac {9} {4} \tau  ) \|Y_{h,m}\|^2 
+ 
2 \tau\sum_{i=1}^m\|\nabla Y_{h,i}\|^2
\\
& \leq
\|Y_{h,0}\|^2
+
\tfrac {9\tau}{4} \sum_{i=1}^{m-1}\|Y_{h,i}\|^2
+
4 \tau\sum_{i=1}^m\|F(W_{A_h}^i)\|^2
\\
& \leq
\|Y_{h,0}\|^2
+
9\tau \sum_{i=1}^{m-1} 
\! \Big ( \!
\tfrac14 \|Y_{h,i}\|^2
+
2 \tau\sum_{j=1}^i \|\nabla Y_{h,j}\|^2
\Big )
+
4 \tau\sum_{i=1}^m\|F(W_{A_h}^i)\|^2.
\end{split}
\end{equation}
}}
By virtue of the Gronwall inequality, we infer that
 \begin{align}
\|Y_{h,m}\|^2
+
\tau\sum_{i=1}^m\|\nabla Y_{h,i}\|^2
\leq
C
\|Y_{h,0}\|^2
+
C\tau\sum_{i=1}^m\|F(W_{A_h}^i)\|^2.
 \end{align}
Let $\mathcal{S}_{\tau,h}(t)=\mathcal{S}_{\tau,h}^i,$ for $t\in[t_{i-1},t_i)$ and by $\chi_B$ we denote the characteristic function
of a set $B\subset \mathbb{R}$. Then $W_{A_h}^m$ can be reformulated as $W_{A_h}^m=\int_0^{T}\chi_{[0, t_m)}(s)
\mathcal{S}_{\tau,h}(t_m-s) P_h \,\dd W(s)$.  As in \eqref{eq:boundness-f(WAH)},  employing \eqref{eq:H1-L6}, \eqref{eq:relation-A-Ah},
\eqref{eq:boundness-full-operator-sum} and Burkholder-Davis-Gundy-type inequality helps us to deduce
{\color{black}{
 \begin{align}\label{eq:regularity-W(Ah)}
 \begin{split}
 \|F(W^m_{A_h})\|_{L^{2p}(\Omega;H)}
 &\leq
 C
 \big (
 1
 +
\big \| W_{A_h}^m \big \|
 ^{3}_{L^{6p}(\Omega;L_6)}
 \big)
 \\
 &\leq
 C
 \big (
 1
 +
\big \| W_{A_h}^m \big \|
 ^{3}_{L^{6p}(\Omega;\dot{H}^{\frac d3})}
 \big)
 \\
&\leq
C\Big(1+
 \Big(\int_0^T\|\chi_{[0, t_m)}(s)A_h^{\frac d6}\mathcal{S}_{\tau,h}(t_m-s)P_hQ^{\frac12}\|_{\mathcal{L}_2}^2\,\dd s\Big)^{3/2}\Big)
 \\
 &
 \leq
 C \Big(
 1 +
 \Big(
 \tau\sum_{i=0}^{m-1}\|A_h^{\frac d6}\mathcal{S}_{\tau,h}^{m-i}P_hQ^{\frac12}\|^2_{\mathcal{L}_2}
 \Big)^{3/2}
 \Big)
 \leq
 C(1+\|A^{\frac{d-3}6}Q^{\frac12}\|_{\mathcal{L}_2}^{3}) <\infty,
 \end{split}
 \end{align}
 }}
 for any $m\in\{1,2,\cdots,M\}$.
 This together with Assumption \ref{assum:initial-value-u0} shows \eqref{eq:boundness-full-FEM}.
$\square$

%

Next we  prove Theorem \ref{them:error-estimate-full-FME}.

{\it Proof of Theorem \ref{them:error-estimate-full-FME}.}
Similarly to the semi-discrete case, by introducing the  auxiliary problem,
\begin{align}\label{eq:auxiliary-problem-full}
\begin{split}
\widetilde{X}_{h,m} - \widetilde{X}_{h,m-1} +  \tau A_h\widetilde{X}_{h,m} =  \tau P_hF ( X(t_m) ) +  P_h \Delta W_m,\; \widetilde{X}_{h,0}=P_hX_0,
\end{split}
\end{align}
whose solution can be recasted as
\begin{align}\label{eq:solution-full-auxiliary problem}
\begin{split}
\widetilde{X}_{h,m}
=
\mathcal{S}^m_{\tau,h}P_hX_0
+
\tau\sum_{i=0}^{m-1}\mathcal{S}^{m-i}_{\tau,h}P_hF(X(t_{i+1}))
+
W_{A_h}^m,
\end{split}
\end{align}
we decompose the considered error term $\| X(t_m)-X_{h,m} \|_{L^{2p}(\Omega;H)}$ into two parts:
\begin{equation}\label{eq:full-error-decomp}
\| X(t_m) - X_{h,m} \|_{L^{2p}(\Omega;H)}
\leq
\| X(t_m) - \widetilde{X}_{h, m} \|_{L^{2p}(\Omega;H)}
+
\| \widetilde{X}_{h, m}  -  X_{h,m} \|_{L^{2p}(\Omega;H)}.
\end{equation}
Resorting to  \eqref{eq:boundeness-f}, \eqref{eq:H1-L6}, \eqref{eq:relation-A-Ah}, \eqref{eq:boundness-full-discrete-operator} and
\eqref{eq:boundness-full-operator-sum},
one can infer that, for any $m\in\{1,2,\cdots, M\}$,
{\color{black}{\begin{align}\label{eq:boundness-full-discrete-auxiliary-problem}
\begin{split}
&
\|\widetilde{X}_{h,m}-W_{A_h}^m\|_{L^{2p}(\Omega;L_6)}
\leq
\|\widetilde{X}_{h,m}-W_{A_h}^m\|_{L^{2p}(\Omega;\dot{H}^{\frac d3})}
\\
\quad &
\leq
\|\mathcal{S}^m_{\tau,h}P_hX_0\|_{L^{2p}(\Omega;\dot{H}^{\frac d3})}
+
\tau\sum_{i=0}^{m-1}\|\mathcal{S}^{m-i}_{\tau,h}P_hF(X(t_{i+1}))\|_{L^{2p}(\Omega;\dot{H}^{\frac d3})}
\\
&
\leq
C\|X_0\|_{L^{2p}(\Omega;\dot{H}^{\frac d3})}
+
C\sup_{s\in[0,T]}\|F(X(s))\|_{L^{2p}(\Omega;H)}\tau\sum_{i=1}^m t_{m-i}^{-\frac d6}
<\infty
,
\end{split}
\end{align}
which together with the fact $\| W_{A_h}^m\|_{L^{2p}(\Omega;L_6)} < \infty$, implied by \eqref{eq:regularity-W(Ah)},  yields
\begin{align}
\|\widetilde{X}_{h,m}\|_{L^{2p}(\Omega;L_6)}<\infty.
\end{align}}}
As the first step, we aim to bound the error $\| X(t_m) - \widetilde{X}_{h, m} \|_{L^{2p}(\Omega;H)} $.
Subtracting \eqref{eq:solution-full-auxiliary problem} from \eqref{eq:intro-mild-solution},
the error $X(t_m)-\widetilde{X}_{h,m}$ can be splitted into the following three  terms:
\begin{align}\label{eq:full-error-decompostion}
\|X(t_m)-\widetilde{X}_{h,m}\|_{L^{2p}(\Omega;H)}
&
=
\|(\mathcal{S}(t_m)-\mathcal{S}_{\tau,h}^mP_h)X_0\|_{L^{2p}(\Omega;H)}
\nonumber \\
&
\quad +
\Big\|\int_0^{t_m}\mathcal{S}(t_m-s)F(X(s))\,\dd s - \tau\sum_{i=0}^{m-1}\mathcal{S}^{m-i}_{\tau,h}P_hF(X(t_{i+1}))\Big\|_{L^{2p}(\Omega;H)}
\nonumber \\
&
\quad
+
\Big\|\int_0^{t_m}\mathcal{S}(t_m-s)\,\dd W(s)- \sum_{i=0}^{m-1}\mathcal{S}^{m-i}_{\tau,h}P_h \Delta W_{i+1} \Big\|_{L^{2p}(\Omega;H)}
\nonumber \\
&
: =
J_1+J_2+J_3.
\end{align}
In the same manner as \eqref{eq:estiamte-I1}, the first term $J_1$ can be estimated with the aid of \eqref{eq:error-estimates-deter-full-intial-positive-norm},
\begin{align}
J_1
\leq
C(h^\gamma + \tau^{\frac\gamma2})\|X_0\|_{L^{2p}(\Omega;\dot{H}^\gamma)}.
\end{align}
To treat the term $J_2$, we decompose it into two terms as follows:
\begin{align}\label{eq:decompose-J2}
\begin{split}
J_2
&
\leq
\Big\|\sum_{i=0}^{m-1}\int_{t_i}^{t_{i+1}}\mathcal{S}(t_m-s)(F(X(s))-F(X(t_{i+1}))\,\dd s\Big\|_{L^{2p}(\Omega;H)}
\\
&
\quad
+
\Big\| \sum_{i=0}^{m-1}\int_{t_i}^{t_{i+1}}(\mathcal{S}(t_m-s)-\mathcal{S}^{m-i}_{\tau,h}P_h)
F(X(t_{i+1}))\,\dd s \Big\|_{L^{2p}(\Omega;H)}
\\
&
: =
J_{21}+J_{22}.
\end{split}
\end{align}
Since the term $J_{22}$ is easy, we treat it first.  Performing standard variable transformations $t_m-s=\sigma$, $m-i=j$
and using \eqref{eq:boundeness-f}, \eqref{prop:boundness-f(x)-f(y)-H}, \eqref{eq:error-estimates-deter-full-intial-positive-norm}
and \eqref{lem:eq-error-full-sum}  yield
\begin{align}\label{eq:boundness-L22}
J_{22}
=&
\big\|\sum_{j=1}^{m}\int_{t_{j-1}}^{t_j}(\mathcal{S}(\sigma)-\mathcal{S}^j_{\tau,h}P_h)F(X(t_{m-j+1}))\,\dd \sigma
\big\|_{L^{2p}(\Omega;H)}
\nonumber
\\
\leq
&
\Big\|\int_0^{t_m}\Psi_{\tau,h}(\sigma)F(X(t_m))\,\dd \sigma\Big\|_{L^{2p}(\Omega;H)}
\nonumber
\\
&
+
\sum_{j=1}^{m}\int_{t_{j-1}}^{t_j}\big\|\Psi_{\tau,h}(\sigma)\left(F(X(t_{m-j+1}))
-
F(X(t_m))\right)\big\|_{L^{2p}(\Omega;H)}\,\dd \sigma
\nonumber
\\
\leq
&
C(h^2+\tau)\|F(X(t_m))\|_{L^{2p}(\Omega;H)}
\nonumber
\\
&
+
C\sum_{j=1}^{m}\int_{t_{j-1}}^{t_j}(h^\gamma+\tau^{\frac\gamma2})\sigma^{-\frac\gamma2}\|F(X(t_{m-j+1}))
-
F(X(t_m))\big\|_{L^{2p}(\Omega;H)}\,\dd \sigma
\nonumber
\\
\leq
&
C(h^2+\tau)\sup_{s\in[0,T]}\|F(X(s))\|_{L^{2p}(\Omega;H)}
+
C\sum_{j=1}^{m}\int_{t_{j-1}}^{t_j}(h^\gamma+\tau^{\frac\gamma2})\sigma^{-\frac\gamma2}t_{j-1}^{\alpha_\gamma}\,\dd \sigma
\nonumber
\\
\leq
&
C(h^\gamma+\tau^{\frac\gamma2}),
\end{align}
where for any fixed number $\delta_0\in(\frac32,2)$, $\alpha_\gamma = 0 $ for $ \gamma \in {\color{black}{[\frac d3, \delta_0]}} $
and $\alpha_\gamma = \tfrac12 $ for $ \gamma \in (\delta_0, 2 ] $ by \eqref{prop:boundness-f(x)-f(y)-H}.
In the next step, we start the estimate of $J_{21}$.
Noting that, for $s\in[t_i,t_{i+1})$
\begin{align}
X(t_{i+1})
=
\mathcal{S}(t_{i+1}-s)X(s)
+
\int_{s}^{t_{i+1}}\mathcal{S}(t_{i+1}-\sigma)F(X(\sigma))\,\dd \sigma
+
\int_{s}^{t_{i+1}}\mathcal{S}(t_{i+1}-\sigma)\,\dd W(\sigma),
\end{align}
and thus  using the Taylor formula {\color{black}{helps}} us to split $J_{21}$ into four terms:
\begin{align}
\begin{split}
J_{21}
\leq
&
\Big\|
\sum_{i=0}^{m-1}\int_{t_i}^{t_{i+1}}\mathcal{S}(t_m-s)F'(X(s))(\mathcal{S}(t_{i+1}-s)-I)X(s)\,\dd s
\Big\|_{L^{2p}(\Omega;H)}
\\
&+
\Big\|
\sum_{i=0}^{m-1}\int_{t_i}^{t_{i+1}}\mathcal{S}(t_m-s)F'(X(s))\int_s^{t_{i+1}}\mathcal{S}(t_{i+1}-\sigma)F(X(\sigma))\,\dd \sigma\,\dd s
\Big\|_{L^{2p}(\Omega;H)}
\\
&+
\Big\|
\sum_{i=0}^{m-1}\int_{t_i}^{t_{i+1}}\mathcal{S}(t_m-s)F'(X(s))\int_s^{t_{i+1}}\mathcal{S}(t_{i+1}-\sigma)\,\dd W(\sigma)\,\dd s
\Big\|_{L^{2p}(\Omega;H)}
\\
&+
\Big\|
\sum_{i=0}^{m-1}\int_{t_i}^{t_{i+1}}\mathcal{S}(t_m-s)R_F(X(s),X(t_{i+1}))\,\dd s
\Big\|_{L^{2p}(\Omega;H)}
\\
: =
&
J_{21}^1+J_{21}^2+J_{21}^3+J_{21}^4.
\end{split}
\end{align}
Here the remainder term $R_F $ reads,
\begin{equation}
\begin{split}
& R_F(X(s),X(t_{i+1}))
\\ &
\quad
:=
\int_0^1
F'' \big( X(s)+\lambda(X(t_{i+1})-X(s)) \big) \big( X(t_{i+1})-X(s), X(t_{i+1})-X(s) \big) (1-\lambda)\,\dd \lambda.
\end{split}
\end{equation}
In the sequel we treat the above four terms one by one.
Thanks to  \eqref{spatial-temporal-S(t)}, \eqref{eq:relation-L(1)-H(-deta))},  \eqref{them:spatial-regularity-mild-stoch}, 
\eqref{eq:thm-wellposed-L6}, \eqref{eq:definition-derivative-F} and H\"{o}lder's inequality,
we derive, for $\gamma\in[\tfrac{d}{3},2]$ and any fixed $\delta_0 \in (\tfrac32, 2 )$,
\begin{align}\label{eq:estimate-L2111}
\begin{split}
J_{21}^1
&\leq
C\sum_{i=0}^{m-1}\int_{t_i}^{t_{i+1}}(t_m-s)^{-\frac{\delta_0}2}
      \|A^{-\frac{\delta_0}2}F'(X(s)) (\mathcal{S}(t_{i+1}-s)-I)X(s)\|_{L^{2p}(\Omega;H)}\,\dd s
\\
&\leq
C\sum_{i=0}^{m-1}\int_{t_i}^{t_{i+1}}(t_m-s)^{-\frac{\delta_0}2}\|F'(X(s))
    (\mathcal{S}(t_{i+1}-s)-I)X(s)\|_{L^{2p}(\Omega;L_1)}\,\dd s
\\
&\leq
C\sum_{i=0}^{m-1}\int_{t_i}^{t_{i+1}}(t_m-s)^{-\frac{\delta_0}2}
      {\color{black}{\big( 1 + \| X (s) \|^2_{ L^{8p}(\Omega; L_4) } \big)}}
      \|(\mathcal{S}(t_{i+1}-s)-I)X(s)\|_{L^{4p}(\Omega;H)}\,\dd s
\\
&\leq
C\tau^{ \frac{\gamma}{2} }\sum_{i=0}^{m-1}\int_{t_i}^{t_{i+1}}(t_m-s)^{-\frac{\delta_0}2}\,\dd s
{\color{black}{\Big( 1 + \sup_{s\in[0,T]} \| X (s) \|^2_{ L^{8p}(\Omega; L_4) } \Big)}}
\sup_{s\in[0,T]}\|X(s)\|_{L^{4p}(\Omega;\dot{H}^\gamma)}
\\
& \leq
 C\tau^{ \frac{\gamma}{2} }.
\end{split}
\end{align}
For the second term $J_{21}^2$, using \eqref{spatial-temporal-S(t)},  \eqref{eq:relation-L(1)-H(-deta))}, \eqref{eq:boundeness-f},
\eqref{eq:thm-wellposed-L6} and  \eqref{eq:definition-derivative-F}
implies, for any fixed $\delta_0 \in (\tfrac32, 2 )$
\begin{align}\label{eq:estimate-L2112}
\begin{split}
J_{21}^2
&\leq
\sum_{i=0}^{m-1}\int_{t_i}^{t_{i+1}}\int_s^{t_{i+1}}(t_m-s)^{-\frac{\delta_0}2}
\big\|A^{-\frac{\delta_0}2}F'(X(s))\mathcal{S}(t_{i+1}-\sigma)F(X(\sigma))\big\|_{L^{2p}(\Omega;H)}\,\dd \sigma\,\dd s
\\
&\leq
\sum_{i=0}^{m-1}\int_{t_i}^{t_{i+1}}\int_{s}^{t_{i+1}}(t_m-s)^{-\frac{\delta_0}2}
\big\|F'(X(s))\mathcal{S}(t_{i+1}-\sigma)F(X(\sigma))\big\|_{L^{2p}(\Omega;L_1)}\,\dd \sigma\,\dd s
\\
&\leq
C\sum_{i=0}^{m-1}\int_{t_i}^{t_{i+1}}\int_s^{t_{i+1}}(t_m-s)^{-\frac{\delta_0}2}
{\color{black}{\big( 1 + \| X (s) \|^2_{ L^{4p}(\Omega; L_4) } \big)}}
\|F(X(\sigma))\|_{L^{4p}(\Omega;H)}\,\dd \sigma\,\dd s
\\
&\leq
C\tau \int_0^{t_m}(t_m-s)^{-\frac{\delta_0}2} \dd s 
{\color{black}{\Big( 1 + \sup_{s\in[0,T]} \| X (s) \|^2_{ L^{4p}(\Omega; L_4) } \Big)}}
\sup_{s\in[0,T]}\|F(X(s))\|_{L^{4p}(\Omega;H)}
\\
& \leq
 C\tau.
\end{split}
\end{align}
To estimate $J_{21}^3$,
we first apply the stochastic Fubini theorem (e.g. see
\cite[Theorem 4.18]{da2014stochastic}) and  the Burkholder-Davis-Gundy-type inequality  to obtain
\begin{align}
\begin{split}
J_{21}^3
=
&
\Big\|\sum_{i=0}^{m-1}\int_{t_i}^{t_{i+1}}\int_{t_i}^{t_{i+1}}\chi_{[s,t_{i+1})}(\sigma)
   \mathcal{S}(t_m-s)F'(X(s))\mathcal{S}(t_{i+1}-\sigma)\,\dd W(\sigma)\,\dd s \Big\|_{L^{2p}(\Omega;H)}
\\
=
&
\Big\|\sum_{i=0}^{m-1}\int_{t_i}^{t_{i+1}}\int_{t_i}^{t_{i+1}}\chi_{[s,t_{i+1})}(\sigma)
   \mathcal{S}(t_m-s)F'(X(s))\mathcal{S}(t_{i+1}-\sigma)\,\dd s \, \dd W(\sigma) \Big\|_{L^{2p}(\Omega;H)}
\\
\leq
&
C
\Big(\sum_{i=0}^{m-1}\int_{t_i}^{t_{i+1}}\Big\|\int_{t_i}^{t_{i+1}}
\mathcal{S}(t_m-s)F'(X(s))\chi_{[s,t_{i+1})}(\sigma)\mathcal{S}(t_{i+1}-\sigma)Q^{\frac12}\,\dd s\Big\|_{L^{2p}
(\Omega;\mathcal{L}_2)}^2\,\dd \sigma\Big)^\frac12.
\end{split}
\end{align}
Further, we employ the H\"{o}lder inequality, \eqref{eq:definition-derivative-F}, \eqref{eq:thm-wellposed-L6},
{\color{black}{the Sobolev embedding inequality
$\dot{H}^{\frac d3}\subset L_6(\mathcal{D}), d \in \{1,2,3 \}$}}  and \eqref{spatial-temporal-integrand-S(t)} 
with ${\color{black}{\rho = \max\{0,\frac d3-\gamma+1\}}}$ to get
\begin{equation}
\begin{split}
J_{21}^3
\leq
&
C\tau^\frac12
\Big(\sum_{i=0}^{m-1}\int_{t_i}^{t_{i+1}}\int_{t_i}^{t_{i+1}}\sum_{j=1}^\infty\|
\mathcal{S}(t_m-s)F'(X(s))\mathcal{S}(t_{i+1}-\sigma)Q^{\frac12}\eta_j\|_{L^{2p}(\Omega;H)}^2\,\dd s\,\dd \sigma\Big)^\frac12
\\
\leq
&
C\tau^\frac12 \Big(\sum_{i=0}^{m-1}\int_{t_i}^{t_{i+1}}
{\color{black}{\big( 1 + \| X (s) \|^4_{ L^{4p}(\Omega; L_6) } \big)}}
\,\dd s
\sum_{j=1}^\infty\int_{t_i}^{t_{i+1}}\|\mathcal{S}(t_{i+1}-\sigma)Q^{\frac12}\eta_j\|_{L_6}^2\,\dd \sigma\Big)^\frac12
\\
\leq
&
C\tau \Big(\sum_{i=0}^{m-1}
{\color{black}{\Big( 1 + \sup_{s\in[0,T]} \| X (s) \|^4_{ L^{4p}(\Omega; L_6) } \Big)}}
\sum_{j=1}^\infty
\int_{t_i}^{t_{i+1}}\|A^{\frac d6}\mathcal{S}(t_{i+1}-\sigma)Q^{\frac12} \eta_j \|^2\,\dd \sigma\Big)^\frac12
\\
\leq
&
C\tau \Big(\sum_{i=0}^{m-1}
\sum_{j=1}^\infty
\int_{t_i}^{t_{i+1}}\|A^{\frac12 ( \frac d 3 - \gamma + 1 ) }\mathcal{S}(t_{i+1}-\sigma) A^{ \frac{\gamma - 1}{2} } Q^{\frac12} \eta_j \|^2\,\dd \sigma\Big)^\frac12
\\
\leq
&
C {\color{black}{\tau^{\frac{2-\max\{0,\frac d 3-\gamma+1\}}2}}}
{\color{black}{\|A^{\frac{\gamma-1}2}Q^{\frac12}\|_{\mathcal{L}_2}}}
\leq
C {\color{black}{\tau^{\frac\gamma2}}},
\end{split}
\end{equation}
where $\eta_j, j \in \N$ is any ON-basis of $H$ and the last inequality holds due to $\gamma \in [ \tfrac d 3, 2 ]$ and $ d \leq 3$.
%
%
%
At the moment we are in a position to bound the term $J_{21}^4$.
Owing to \eqref{spatial-temporal-S(t)} with $\nu=\delta_0\in (\frac32,2)$ and using \eqref{eq:definition-derivative-F}, 
\eqref{them:temporal-regularity-mild-stoch}, \eqref{eq:thm-wellposed-L6},
\eqref{eq:H1-L6} and H\"{o}lder's inequality, we learn that
\begin{align}
\begin{split}
J_{21}^4
&\leq
C\sum_{i=0}^{m-1}\int_{t_i}^{t_{i+1}}(t_m-s)^{-\frac{\delta_0}2}
\|A^{-\frac{\delta_0}2}R_F(X(s),X(t_{i+1}))\|_{L^{2p}(\Omega;H)}\,\dd s
\\
&
\leq
C\sum_{i=0}^{m-1}\int_{t_i}^{t_{i+1}}(t_m-s)^{-\frac{\delta_0}2}
     \|R_F(X(s),X(t_{i+1}))\|_{L^{2p}(\Omega;L_1)}\,\dd s
\\
&
\leq
C\sum_{i=0}^{m-1}\int_{t_i}^{t_{i+1}}(t_m-s)^{-\frac{\delta_0}2}
\int_0^1
\big\|\;\|X(t_{i+1})-X(s)\|
\\
&
\qquad
\,
\times 
{\color{black}{
\big[
( 1-\lambda ) \| X(s) \|_{L_4} + \lambda  \| X(t_{i+1}) \|_{L_4}
\big]
}}
\,\|X(t_{i+1})-X(s)\|_{L_4}\;\big\|_{L^{2p}(\Omega;\mathbb{R})}
\, \dd \lambda
\,\dd s
\\
&
\leq
C\sum_{i=0}^{m-1}
       \int_{t_i}^{t_{i+1}}(t_m-s)^{-\frac{\delta_0}2}\|X(t_{i+1})-X(s)\|_{L^{8p}(\Omega;H)}\|X(t_{i+1})-X(s)\|_{L^{8p}
       (\Omega;{\color{black}{\dot{H}^{d/3}}})}\,\dd s
\\
&\qquad
\times
{\color{black}{\sup_{s\in[0,T]} \| X ( s ) \|_{ L^{ 4 p} ( \Omega; L_4 ) } }}
\\
&\leq
C
{\color{black}{ \tau^{\frac{\min\{1,\gamma\}+\min\{1,\gamma-\frac d3\}}2}
}}
         \sum_{i=0}^{m-1}\int_{t_i}^{t_{i+1}}(t_m-s)^{-\frac{\delta_0}2}\,\dd s
\\
& \leq
C \tau^{\frac\gamma2}.
\end{split}
\end{align}
Putting the above four estimates together results in
\begin{align}
J_{21}\leq C\tau^{\frac\gamma2},
\end{align}
which together with \eqref{eq:decompose-J2} and \eqref{eq:boundness-L22} shows
\begin{align}
J_2
\leq
C(h^\gamma+\tau^{\frac\gamma2}).
\end{align}
Concerning the term $J_3$,  \eqref{lem:eq-error-determ-full-sum-square-norm},    \eqref{eq:A-Q-condition}
and the Burkholder-Davis-Gundy type inequality  show
\begin{align}
\begin{split}
J_3
=
&
\Big\|
\int_0^{t_m}\Psi_{\tau,h}(t_m-s)\,\dd W(s)
\Big\|_{L^{2p}(\Omega;H)}
\leq
C_p
\Big(
\int_0^{t_m}\|\Psi_{\tau,h}(t_m-s)Q^{\frac12}\|_{\mathcal{L}_2}^2\,\dd s
\Big)^{\frac12}
\\
=
&
C_p
\Big(
\int_0^{t_m}\|\Psi_{\tau,h}(s)Q^{\frac12}\|_{\mathcal{L}_2}^2\,\dd s
\Big)^{\frac12}
\leq
C(h^\gamma+\tau^{\frac\gamma2})\|A^{\frac{\gamma-1}2}Q^{\frac12}\|_{\mathcal{L}_2}.
\end{split}
\end{align}
Gathering the above three estimates together implies
\begin{align}\label{eq:error-between-true-auxiliary-problem}
\|X(t_m)-\widetilde{X}_{h,m}\|_{L^{2p}(\Omega;H)}
\leq
C(h^\gamma+\tau^{\frac\gamma2}).
\end{align}

Next we turn our attention to the estimate of $\widetilde{e}_m:=\widetilde{X}_{h,m}-X_{h,m}$, which obeys
\begin{align}\label{eq:em-problem}
\begin{split}
\frac{\widetilde{e}_m-\widetilde{e}_{m-1}}\tau + A_h\widetilde{e}_m= P_h(F(X(t_m))-F(X_{h,m})),\; \widetilde{e}_0=0.
\end{split}
\end{align}
By multiplying this equation by $\widetilde{e}_m$,
one can observe
\begin{align}\label{eq:error-equation-controll}
\begin{split}
&
\tfrac12(\|\widetilde{e}_m\|^2-\|\widetilde{e}_{m-1}\|^2)
+
\tau\left<\nabla \widetilde{e}_m,\nabla \widetilde{e}_m\right>
\\
& \quad
\leq
\tau \big \langle F(\widetilde{X}_{h,m})-F(X_{h,m}),\widetilde{e}_m \big \rangle
+
\tau \big \langle F(X(t_m))-F(\widetilde{X}_{h,m}),\widetilde{e}_m \big \rangle.
\end{split}
\end{align}
Here we also used the definition of $A_h$ in  \eqref{eq:definition-discrete-A} and  the fact
$\frac12(\|\widetilde{e}_m\|^2-\|\widetilde{e}_{m-1}\|^2)\leq \big<\widetilde{e}_m-\widetilde{e}_{m-1},\widetilde{e}_m\big>$.
Thanks to \eqref{eq:local-lipschitz-condition-F} and  \eqref{eq:relation-L(65)-H(-1)},
\begin{align}
\begin{split}
\tfrac12(\|\widetilde{e}_m\|^2&-\|\widetilde{e}_{m-1}\|^2)
+
\tau\left<\nabla \widetilde{e}_m,\nabla \widetilde{e}_m\right>
\\
&\leq
\tau\|\widetilde{e}_m\|^2
+\tau\|A^{-\frac12} (F(X(t_m))-F(\widetilde{X}_{h,m}))\|\,\|\nabla \widetilde{e}_m\|
\\
&\leq
\tau\|\widetilde{e}_m\|^2
+\tfrac\tau2\|F(X(t_m))-F(\widetilde{X}_{h,m})\|_{L_{\frac65}}^2
+
\tfrac\tau2\|\nabla \widetilde{e}_m\|^2
\\
&\leq
\tau\|\widetilde{e}_m\|^2
+
C\tau\|X(t_m)-\widetilde{X}_{h,m}\|^2(1+\|X(t_m)\|_{L_6}^4+\|\widetilde{X}_{h,m}\|_{L_6}^4)
+
\tfrac\tau2\|\nabla \widetilde{e}^m\|^2.
\end{split}
\end{align}
Since Lemma \ref{lem:boundess-full-problem} and \eqref{eq:boundness-full-discrete-auxiliary-problem} ensure
$
\mathbf{E} [ \|\widetilde{e}_m\|^{2p} ]
+
\tau\sum_{i=1}^m \mathbf{E} [\|\nabla \widetilde{e}_i\|^2 ]
<
\infty,
$
by summation on $m$ and  calling the Gronwall inequality and the fact  $\widetilde{e}_0=0$, it holds
\begin{align}
\begin{split}
\|\widetilde{e}_m\|^2
\leq
C\tau\sum_{i=1}^m\|X(t_i)-\widetilde{X}_{h,i}\|^2(1+\|X(t_i)\|^4_{L_6}+\|\widetilde{X}_{h,i}\|^4_{L_6}).
\end{split}
\end{align}
Therefore,
\begin{align}\label{eq:error-discrete-and-auxiliary-control-by-auxiliary-error-full}
\begin{split}
\|\widetilde{e}_m\|_{L^{2p}(\Omega;H)}
\leq
&
C \tau \sum_{i=1}^m
\|X(t_i)-\widetilde{X}_{h,i}\|_{L^{4p}(\Omega;H)}
\big(
1 +  \|X( t_i )\|^2_{L^{8p}(\Omega;L_6)}
+
\|\widetilde{X}_{h,i}\|^2_{L^{8p}(\Omega;L_6)}
\big)
\\
\leq
&
C(h^\gamma+ \tau^{\frac\gamma2}),
\end{split}
\end{align}
which together with \eqref{eq:error-between-true-auxiliary-problem} shows \eqref{them:error-result-full-FME}
and thus finishes the proof. $\square$
\section{Numerical experiments}
\label{sect:numerical-examples}
In this section, some numerical examples are included to  {\color{black}{illustrate}} the previous findings.
To this end, we consider the following stochastic Allen-Cahn equation in one space dimension
\begin{align}\label{eq:example-SACE}
\left\{\begin{array}{ll}
\frac{\partial u}{\partial t}
=
 \frac{\partial^2 u}{\partial x^2}
 +
 u-u^3
 +
 \dot{W}
 \,&t\in(0,1],\;x\in(0,1),
\\
u(0,x)=sin(\pi x),
&
x\in(0,1),
\\
u(t,0)=u(t,1) = 0, &
t\in(0,1].
\end{array}
\right.
\end{align}
Here $\{ W (t) \}_{t \in[0, 1]}$ stands for a standard $Q$-Wiener process, with two simple choices of covariance operators
$Q = A^{-s}, s \in \{0.5005, 1.5005\}$. 
One can easily see that  Assumption \ref{assum:eq-noise} is fulfilled with $\gamma = 1$ for $Q = A^{-0.5005}$
and $\gamma = 2$ for $Q = A^{-1.5005}$. According to Theorem \ref{them-error-semi-discrete-problem} and
Theorem  \ref{them:error-estimate-full-FME}, the mean-square (MS, $p = 1$) convergence rate in space reads $O( h^{\gamma} )$
and the rate in time $O( \tau^{\frac\gamma2} )$ for $\gamma \in \{ 1, 2 \}$. Since the exact solution is not available,
we turn to fine numerical approximations for reference, using very small  step-sizes $h_{exact}$ and $\tau_{exact}$.
Also, error bounds are always measured in terms of mean-square discretization errors at the endpoint $T = 1$ and
the expectations are approximated by computing averages over 500 samples.

\begin{figure}[!ht]
\centering
      \includegraphics[width=2.6in,height=2.8in] {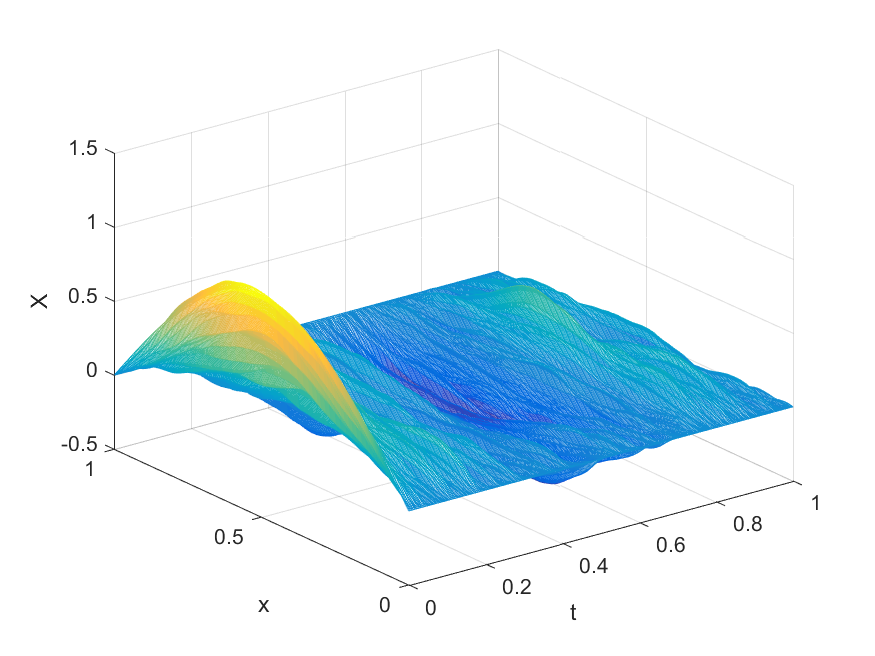}
       \includegraphics[width=2.6in,height=2.8in] {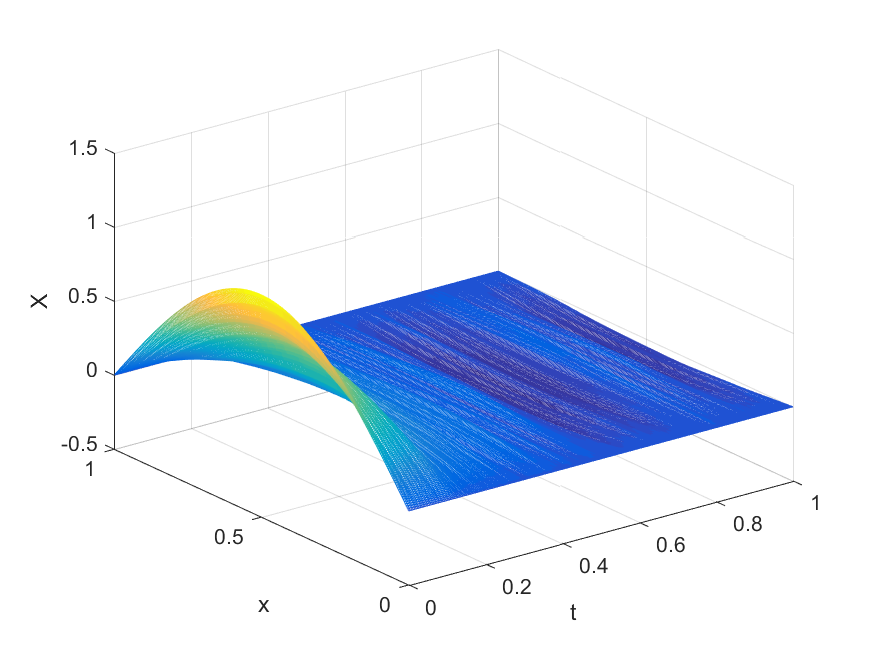}
 \caption{One-sample simulation (Left: $Q=A^{-0.5005}$; Right: $Q=A^{-1.5005}$)  }
\label{fig:path-sample-solution}
\end{figure}

\begin{figure}[!ht]
\centering
      \includegraphics[width=2.6in,height=2.6in] {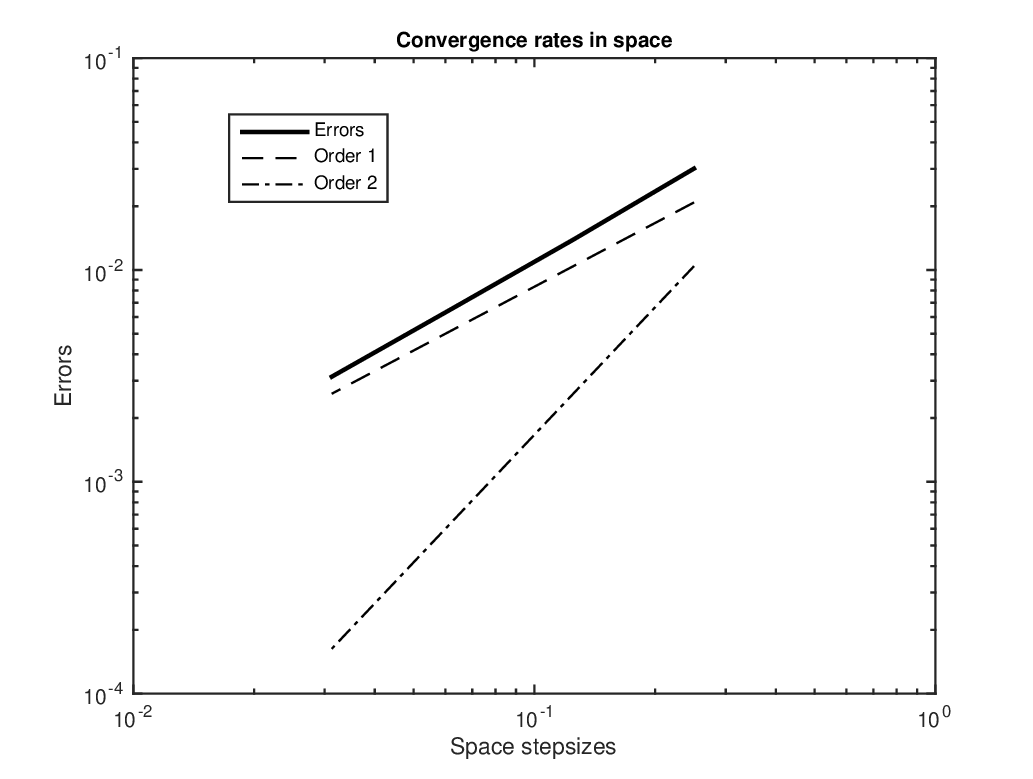}
       \includegraphics[width=2.6in,height=2.6in] {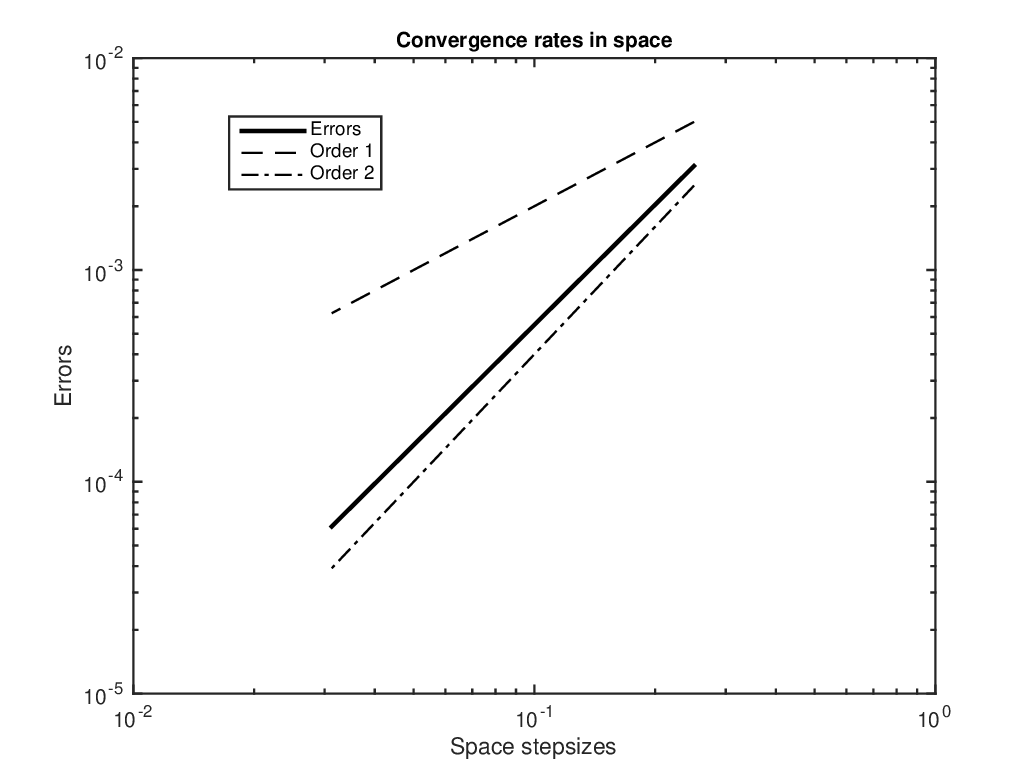}
 \caption{MS convergence rates for spatial discretizations (Left: $Q=A^{-0.5005}$; right: $Q=A^{-1.5005}$)}\label{fig:space-convergence-rate}
\end{figure}

\begin{figure}[!ht]
\centering
      \includegraphics[width=2.6in,height=2.6in] {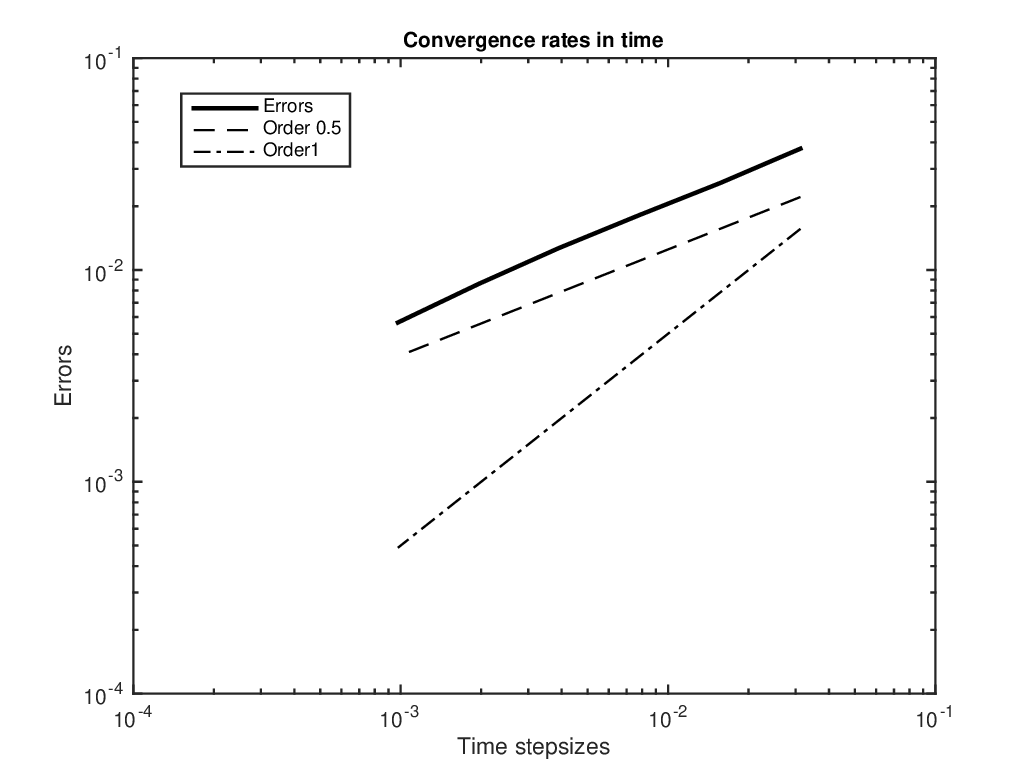}
       \includegraphics[width=2.6in,height=2.6in] {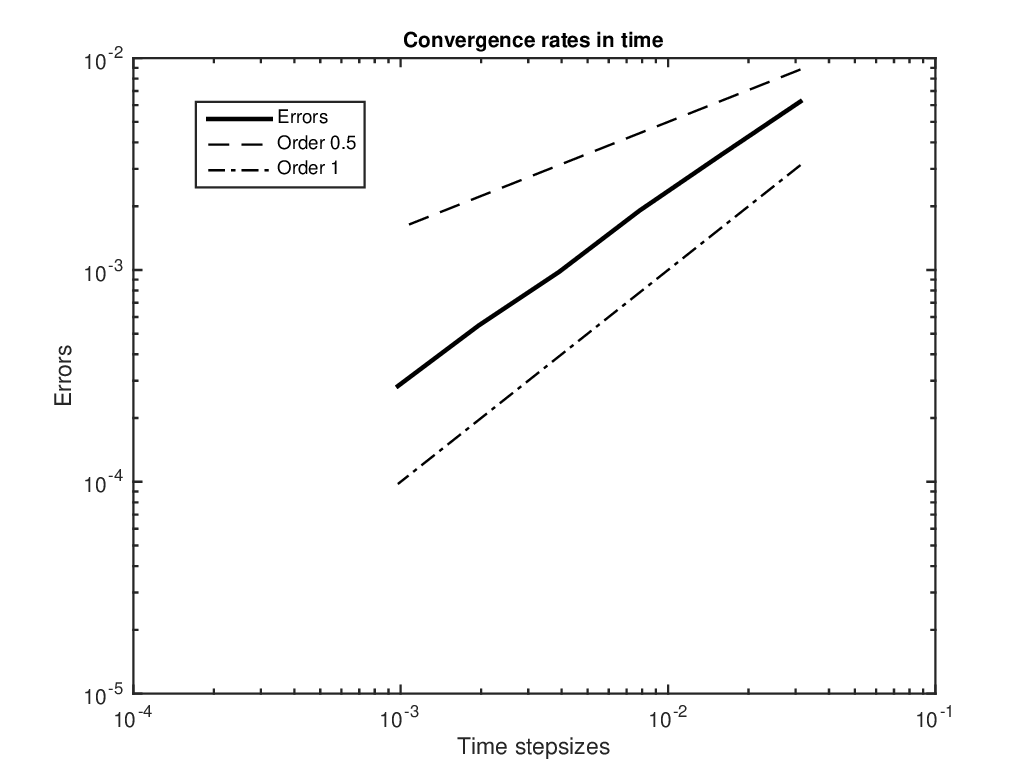}
 \caption{MS convergence rates for time discretizations (Left: $Q=A^{-0.5005}$; right: $Q=A^{-1.5005}$)}
 \label{fig:temporal-convergence-rate}
\end{figure}

In Figure \ref{fig:path-sample-solution}, one-path simulations with $h=\tau=2^{-8}$ are plotted.
There one can observe that the numerical solution behaves more smoothly as the noise becomes  {\color{black}{smoother}}.
To test the convergence rate in space, we perform numerical simulations with four different space step-sizes
$h=2^{-i},  i \in \{2,3,4,5\}$.  The "true solutions"  are computed using $h_{exact}=2^{-7}$, $\tau_{exact}=2^{-15}$.
In Figure \ref{fig:space-convergence-rate}, we depict the spatial errors against space step-sizes and
one can detect the expected convergence rates in space, i.e., order $1$ for $Q = A^{-0.5005}$
and order $2$ for $Q = A^{-1.5005}$.
Lastly, we test the convergence rate in time and take {\color{black}{$h_{exact}=2^{-8}$ and $\tau_{exact}=2^{-14}$.}}
Similarly,  we do numerical approximations with six different time step-sizes {\color{black}{$\tau = 2^{-j}, j \in \{ 5,6,7,8,9,10 \}$}} and
present the resulting errors in Figure \ref{fig:temporal-convergence-rate}.
Clearly,  temporal approximation errors decrease at a slope close to $\frac12$ and $1$ for the above two kinds of noises.
This is consistent with previous theoretical results.

\bibliography{bibfile}
\bibliographystyle{abbrv}

\end{document}